\documentclass[11pt]{amsart}
\usepackage{graphicx}
\usepackage[all]{xy}
\usepackage{amssymb,semantic,a4wide}

\title{Towards a homotopy theory of process algebra}
\author{Philippe Gaucher}
\urladdr{http://www.pps.jussieu.fr/~gaucher/} 
\address{Preuves Programmes et Syst{\`e}mes,
Universit{\'e} Paris 7--Denis Diderot,
Site Chevaleret,
Case 7012,
75205 PARIS Cedex 13,
France}
\thanks{This work has been supported by the ANR ``Invariants
  alg\'ebriques des syst\`emes informatiques'' ANR-05-BLAN-0267.}

\newcommand{\de}{\partial}
\newcommand{\p}\times
\renewcommand{\vec}{\overrightarrow}
\renewcommand{\P}{\mathbb{P}}

\newcommand{\beas}{\begin{eqnarray*}}
\newcommand{\eeas}{\end{eqnarray*}}

\swapnumbers

\newtheorem{theorem}{Theorem}[section]
\newtheorem{proposition}[theorem]{Proposition}
\newtheorem{lemma}[theorem]{Lemma}
\newtheorem{corollary}[theorem]{Corollary}
\newtheorem{definition}[theorem]{Definition}
\newtheorem{nota}[theorem]{Notation}

\newcommand{\bd}{\begin{definition}}
\newcommand{\ed}{\end{definition}}
\newcommand{\bp}{\begin{proposition}}
\newcommand{\ep}{\end{proposition}}
\newcommand{\bth}{\begin{theorem}}
\renewcommand{\eth}{\end{theorem}}
\newcommand{\bpf}{\begin{proof}}
\newcommand{\epf}{\end{proof}}

\newcommand{\fl}[1]{\ar@{->}[ll]_-{#1}}
\newcommand{\fr}[1]{\ar@{->}[rr]^-{#1}}
\newcommand{\fd}[1]{\ar@{->}[dd]_-{#1}}
\newcommand{\fu}[1]{\ar@{->}[uu]^-{#1}}
\newcommand{\f}[2]{\ar@{->}[#1]|{#2}}
\newcommand{\ff}[2]{\ar@2{->}[#1]|{#2}}
\newcommand{\frr}[1]{\ar@{->}[rrrr]^-{#1}}

\renewcommand{\top}{{\mathbf{Top}}}

\newcommand{\cat}{{\mathbf{Cat}}}
\newcommand{\iso}{\cong}
\newcommand{\lp}{\left(}
\newcommand{\rp}{\right)}
\newcommand{\ot}{\otimes}
\newcommand{\vI}{\vec{I}}
\renewcommand{\leq}{\leqslant}
\renewcommand{\geq}{\geqslant}
\newcommand{\brm}[1]{\rm{\mathbf{#1}}}
\newcommand{\dtop}{{\brm{Flow}}}
\newcommand{\set}{{\brm{Set}}}
\newcommand{\poset}{{\brm{PoSet}}}
\newcommand{\proc}{{\brm{Proc}}}
\newcommand{\glob}{{\rm{Glob}}}
\DeclareMathOperator{\rec}{rec}
\DeclareMathOperator{\trans}{Trans}
\DeclareMathOperator{\id}{Id}

\DeclareMathOperator{\cosk}{cosk}
\DeclareMathOperator{\COSK}{\vec{\cosk}}
\newcommand{\liminj}{\varinjlim}

\def\cartesien{%
  \ar@{-}[]+R+<6pt,-2pt>;[]+RD+<6pt,-6pt>%
  \ar@{-}[]+D+<2pt,-6pt>;[]+RD+<6pt,-6pt>%
}
\def\cocartesien{%
  \ar@{-}[]+L+<-6pt,+2pt>;[]+LU+<-6pt,+6pt>%
  \ar@{-}[]+U+<-2pt,+6pt>;[]+LU+<-6pt,+6pt>%
}
\def\hocartesien{%
  \ar@{-}[]+R+<6pt,-2pt>;[]+RD+<6pt,-6pt>_{h}%
  \ar@{-}[]+D+<2pt,-6pt>;[]+RD+<6pt,-6pt>%
}
\def\hococartesien{%
  \ar@{-}[]+L+<-6pt,+2pt>;[]+LU+<-6pt,+6pt>_{h}%
  \ar@{-}[]+U+<-2pt,+6pt>;[]+LU+<-6pt,+6pt>%
}

\makeatletter

\def\varholim@#1#2{%
  \vtop{\m@th\ialign{##\cr
    \hfil$#1\operator@font holim$\hfil\cr
    \noalign{\nointerlineskip\kern1.5\ex@}#2\cr
    \noalign{\nointerlineskip\kern-\ex@}\cr}}%
}
\def\holimproj{%
  \mathop{\mathpalette\varholim@{\leftarrowfill@\textstyle}}\nmlimits@
}
\def\holiminj{%
  \mathop{\mathpalette\varholim@{\rightarrowfill@\textstyle}}\nmlimits@
}

\makeatother

\begin{document}

\begin{abstract} 
  This paper proves that labelled flows are expressive enough to
  contain all process algebras which are a standard model for
  concurrency.  More precisely, we construct the space of execution
  paths and of higher dimensional homotopies between them for every
  process name of every process algebra with any synchronization
  algebra using a notion of labelled flow.  This interpretation of
  process algebra satisfies the paradigm of higher dimensional
  automata (HDA): one non-degenerate full $n$-dimensional cube (no
  more no less) in the underlying space of the time flow corresponding
  to the concurrent execution of $n$ actions.  This result will enable
  us in future papers to develop a homotopical approach of process
  algebras.  Indeed, several homological constructions related to the
  causal structure of time flow are possible only in the framework of
  flows.
\end{abstract}

\maketitle

\tableofcontents

\section{Introduction}

\subsection{Presentation of the results} Process algebras are a
standard way of modelling concurrent processes~\cite{0683.68008,
  MR1365754}.  Some homotopical tools are introduced in~\cite{model3}
to study these concurrent systems.  However the expressiveness of
these tools has never been verified so far. The goal of this paper is
to prove that the category of flows is expressive enough to contain
all process algebras with any synchronization algebra.  This semantics
of process algebra satisfies the paradigm of higher dimensional
automata (HDA)~\cite{fahrenberg05-hda-long, fahrenberg05-hda,
  labelled, Pratt, rvg, exHDA}: one non-degenerate full
$n$-dimensional cube (no more no less) in the underlying space of the
time flow corresponding to the concurrent execution of $n$ actions. In
particular, this semantics of true concurrency does make the
difference between the process $a.b.nil+b.a.nil$ (Figure~\ref{seq})
corresponding to the non-deterministic choice between the sequential
execution of $a$ and $b$, and the sequential execution of $b$ and $a$,
and the process $a.nil|| b.nil$ (Figure~\ref{conc}) corresponding to
the concurrent execution of $a$ and $b$. Figure~\ref{seq} will be
modelled by the boundary of a $2$-cube. Figure~\ref{conc} will be
modelled by a full $2$-cube.

This construction will enable us in future papers to study the
homotopy theory of process algebra, and also the behaviour of new
homological theories introduced in the framework of flows as the
branching and merging homology theories~\cite{exbranch}. These
homology theories are of interest in computer science since they
detect the non-deterministic areas of branchings and mergings of the
time flow of a concurrent process.  Therefore they contain useful
information about the causal structure of the time flow. In
particular, it is likely that they are related in some way to
bisimulation of flows which will be introduced in a future paper too.
It is worth noting that none of the other topological models of
concurrency (local po-spaces~\cite{MR1683333},
d-spaces~\cite{MR2030049}, etc. are convenient for the study of these
homology theories. Indeed, these latter categories contain too many
morphisms, making impossible the construction of \emph{functorial}
branching and merging homology theories: see~\cite[\S 20]{model3}
or~\cite[\S 6]{diCW} for further explanations.  The category of
precubical sets is also not convenient for the study of these homology
theories because of the absence of degenerate cubes (i.e.\ of ``thin''
cubes, that is without volume) and of composition of cubes: see the
introduction and especially Figure~3 of~\cite{Coin} for further
explanations.  However, the construction of this paper does use as an
intermediate category the category of precubical sets.

\subsection{Outline of the paper} Section~\ref{completedescription}
gives a very short presentation of \emph{process algebra} for
mathematicians.  The syntax of the language, as well as the usual
operational semantics (of dimension $1$!) are described.  We take a
version of process algebra without message passing for simplicity.
Section~\ref{slabel} gives the definition of a \emph{decorated
  $\sigma$-labelled precubical set}, where $\sigma$ is a
synchronization algebra.  Our definition of a labelled precubical set
is similar to Goubault's definition~\cite{labelled}, however with some
subtle differences: the new definition takes into account the
synchronization algebra $\sigma$ and the set of labels is not ordered
anymore (see Proposition~\ref{cubeetiquette} and
Notation~\ref{cubeetiquette1}).  This section also presents the
\emph{$\sigma$-labelled directed coskeleton construction}.  The idea
is borrowed from Worytkiewicz's ideas about coskeletal
synchronization~\cite{exHDA}. The reader must know that what we call
\emph{directed coskeleton} is something distinct from the usual
coskeleton. It is defined as a labelled precubical subset of a
labelled coskeleton, the latter coinciding with the usual coskeleton
if the set of labels is a singleton consisting of one action which may
run asynchronously.  Section~\ref{semccsprecube} describes the
denotational semantics of process algebra in terms of
$\sigma$-labelled precubical sets. The only new interpretation is the
one of the parallel composition with synchronization. Several
elementary properties of the parallel composition with synchronization
are then explicitly proved.  Section~\ref{c0} explains the link
between this semantics and the $1$-dimensional operational semantics
of Table~\ref{opsem} and gives an explicit statement corresponding to
the paradigm of higher dimensional automata.  Section~\ref{sflow}
gives the definition of \emph{decorated $\sigma$-labelled flow}, with
a lot of elementary examples.  Section~\ref{rea} is the mathematical
core of the paper. It constructs the \emph{geometric realization
  functor} from precubical sets to flows. The weak S-homotopy model
category introduced in~\cite{model3} is required for this construction
(cf.  Theorem~\ref{pbpb}) and also for the description of some
elementary properties. In particular, it is proved that the geometric
realization of the boundary of the $n$-cube contains a hole with the
correct dimension, that is $n$ (Corollary~\ref{correct}).  Finally,
Section~\ref{reaet} shows how one can associate a $\sigma$-labelled
flow with a $\sigma$-labelled precubical set.
\subsection{Prerequisites} Section~\ref{completedescription} is a
quick introduction about \emph{process algebra} which hopefully
contains enough details. A possible reference is~\cite{MR1365754}.
Model categorical techniques are used in Section~\ref{rea} to
construct the geometric realization functor and to prove some basic
facts about it. The proofs of the theorems in Section~\ref{rea} can be
skipped without problem in a first reading, especially the proof of
Theorem~\ref{nopb} which makes heavy use of homotopical material
coming from~\cite{4eme}.  Possible references for \emph{model
  categories} are~\cite{MR1361887, MR99h:55031},
and~\cite{ref_model2}.  The original reference is~\cite{MR36:6480} but
Quillen's axiomatization is not used in this paper. The axiomatization
from Hovey's book is preferred.

\section{Process algebra}
\label{completedescription}

\subsection{Synchronization algebra} 

Let $\Sigma$ be a non-empty set. Its elements are called
\emph{labels}, \emph{actions}, or \emph{events}.  A
\emph{synchronization algebra} on $\Sigma$ (not containing the
distinct elements $\{0,\bot\}$) consists of a binary commutative
associative operation $\sigma(-,-)$ on $\Sigma\cup\{0,\bot\}$ such
that
\begin{itemize} 
\item $\sigma(a,\bot) = \bot$ for every $a\in \Sigma\cup\{\bot\}$,
\item $\sigma(a,b)=0 \iff a=b=0$, 
\item $\forall a\in \Sigma, \sigma(a,0) = a \hbox{ or }\sigma(a,0) =
  \bot$.  Note necessarily, that one has $\sigma(\bot,0)=\bot$.
\end{itemize}
The label $0$ represents the idle action. The role of $\bot$ is to
specify which pair of actions may synchronize. The equality
$\sigma(a,b) = \bot$ means that $a$ and $b$ cannot synchronize. The
equality $\sigma(a,b) = c\in \Sigma$ means that $a$ and $b$ may
synchronize giving an action relabelled by $\sigma(a,b) = c$. The
equality $\sigma(a,0) = a$ means that $a$ may run asynchronously.  The
equality $\sigma(a,0) = \bot$ means that $a$ cannot run
asynchronously.

\bd The \emph{trivial synchronization algebra} $\sigma$ is the
synchronization algebra satisfying $\sigma(a,b)=\bot$ for all
$(a,b)\in \Sigma$ and $\sigma(a,0) = a$ for all $a\in
\Sigma\cup\{0,\bot\}$. It is denoted by $\bot$. \ed

In pure CCS~\cite{0683.68008}, the set $\Sigma$ contains a distinct
action $\tau$ and the complementary $\Sigma\backslash\{\tau\}$ is
equipped with an involution $a\mapsto \overline{a}$. The
synchronization algebra is defined by the following table:
\begin{center}
\begin{tabular}{c|ccccccc}
$\sigma(-,-)$ & 0 & $a$ & $\overline{a}$ & $b$ &   $\overline{b}$ & \dots & $\tau$ \\
\hline
0 & 0 & $a$ & $\overline{a}$ & $b$ & $\overline{b}$ & \dots & $\tau$ \\
$a$ &$a$& $\bot$ & $\tau$  &  $\bot$ & $\bot$  &  &$\bot$\\
$\overline{a}$ &$\overline{a}$& $\tau$  & $\bot$  & $\bot$  &$\bot$   &  &  $\bot$ \\
$b$ &$b$& $\bot$ & $\bot$  & $\bot$ & $\tau$ &  &  $\bot$  \\
$\overline{b}$ &$\overline{b}$& $\bot$  &  $\bot$  &  $\tau$ &  $\bot$  &  & $\bot$ \\
$\dots$ &$\dots$&&&&&&\\
$\tau$ &$\tau$& $\bot$   & $\bot$   & $\bot$ & $\bot$   &  &  $\bot$.
\end{tabular}
\end{center}
In pure CCS, each pair of actions $a,\overline{a}$ for $a\in
\Sigma\backslash\{\tau\}$ may synchronize to form a synchronization
action labelled by $\tau$, and actions labelled by $\tau$ cannot
synchronize further. All labelled events may occur asynchronously.

In TCSP~\cite{0628.68025}, the set $\Sigma$ contains a distinct action
$\tau$ and the synchronization algebra is defined by the following
table:
\begin{center}
\begin{tabular}{c|ccccccc}
$\sigma(-,-)$ & $0$ &  $a$ & $b$ & \dots & $\tau$ \\
\hline
$0$ & $0$ & $\bot$ & $\bot$ & \dots & $\tau$ \\ 
$a$ & $\bot$ & $a$ & $\bot$  &  $\bot$ & $\bot$ \\
$b$ & $\bot$ & $\bot$ & $b$  & $\bot$ & $\bot$  \\
$\dots$ &&&&&\\
$\tau$ & $\tau$ & $\bot$   & $\bot$   & $\bot$ & $\bot$.  
\end{tabular}
\end{center}
In TCSP, an action $a\in \Sigma\backslash\{\tau\}$ must synchronize
with another action labelled by $a$ to occur in a parallel
composition. So non-$\tau$-labelled events cannot occur
asynchronously.

\subsection{Syntax of the language}

The \emph{process names} are generated by the following syntax:
\[
P::= nil \ |\ a.P \ |\ (\nu a)P \ |\ P + P \ |\ P\|  P \ |\
\rec(x)P(x),\] where $P(x)$ means a process name with one free variable
$x$. The variable $x$ must be \emph{guarded}; that is, it must lie in
a prefix term $a.x$ for some $a\in\Sigma$. The set of process names is
denoted by $\proc_\Sigma$.  The names $nil$, $a.b.nil || 
(c.nil+d.nil)$, $\rec(x) (a.x|| b.nil)$ with $a,b,c,d\in \Sigma$ are
examples of elements of $\proc_\Sigma$.

The process $nil$ corresponds to the idle process. The process $a.P$
corresponds to the sequential execution of $a$ and $P$.  The process
$P+Q$ corresponds to the non-deterministic choice of executing $P$ or
$Q$. The process $P|| Q$ corresponds to the concurrent execution of
$P$ and $Q$, with all possible synchronizations of an action $a$ of
$P$, with an action $b$ of $Q$ if $\sigma(a,b)\neq \bot$. The process
$(\nu a)P$ corresponds to the restriction to a local use of $a$: all
transitions of $P$ labelled with $a$ or with an event $b$
synchronizing with $a$ are removed from $(\nu a)P$. For example, in
pure CCS, $(\nu a)P$ is obtained from $P$ by removing all transitions
labelled by $a$ and $\overline{a}$. Last but not least, the process
$\rec(x)P(x)$ corresponds to the recursive execution of $P(x)$.

\subsection{Operational semantics using labelled transition system}
\label{opsemop}

The following definition is standard:

\bd A \emph{labelled transition system} consists in a set of states
$S$, with initial state $i$, a set of labels $\Sigma$ and a transition
relation $\trans \subset S\p \Sigma\p S$.  \ed

All labelled transition systems are loopless in this paper; that is, if
$(a,u,b)\in\trans$, then $a\neq b$.

\bd A labelled transition system $(S,i,\Sigma,\trans)$
\emph{decorated by process names} is a labelled transition system
together with a set map $d\colon S\rightarrow \proc_\Sigma$ called the
\emph{decoration}.  \ed

If $(S,i,\Sigma,\trans,d)$ is a labelled transition system decorated
by process names, then an element $(a,u,b)$ of $\trans$ is denoted by
$d(a) \stackrel{u} \rightarrow d(b)$.  A labelled decorated transition
system will be identified in Section~\ref{c0} with a $\sigma$-labelled
precubical set of dimension $1$ decorated by process names and with a
distinct initial state (see Theorem~\ref{restrict1}).  Intuitively,
the notation $d(a) \stackrel{u} \rightarrow d(b)$ means that
\emph{$d(a)$ behaves like $d(b)$ after the execution of $u$}.  For
example, the transition $\mu.P \stackrel{\mu} \longrightarrow P$ (see
Table~\ref{opsem}) means that the process $\mu.P$ behaves like $P$
after the execution of $\mu$.

The operational semantics of our language is defined by the rules of
Table~\ref{opsem}~\cite{MR1365754}, with $\mu\in\Sigma$.

\begin{table}
\begin{center}
\fbox{
\begin{tabular}{lr}
 \inference[\bf{Act}]{}{\mu.P \stackrel{\mu}\longrightarrow P}&\\
\inference[\bf{Res}]{P \stackrel{\mu} \longrightarrow P' \ \ \mu\neq a\hbox{ and }\sigma(\mu,a)=\bot}{(\nu a)P \stackrel{\mu} \longrightarrow (\nu a)P'}& \\
\inference[\bf{Sum1}]{P \stackrel{\mu}\longrightarrow P'}{P+Q
  \stackrel{\mu}\longrightarrow P'}&
\inference[\bf{Sum2}]{Q
  \stackrel{\mu}\longrightarrow Q'}{P+Q \stackrel{\mu}\longrightarrow
  Q'} \\
\inference[\bf{Par1}]{P \stackrel{\mu}\longrightarrow P'}{P||  Q
  \stackrel{\mu}\longrightarrow P'|| Q} &
\inference[\bf{Par2}]{Q
  \stackrel{\mu}\longrightarrow Q'}{P||  Q \stackrel{\mu}\longrightarrow
  P||  Q'} \\
\inference[\bf{Com}]{P
  \stackrel{a}\longrightarrow P', Q
  \stackrel{b}\longrightarrow Q',\sigma(a,b)\neq\bot}{P||  Q \stackrel{\sigma(a,b)}\longrightarrow
  P'||  Q'}&\\
\inference[\bf{Rec}]{P(\rec(x)P(x)) \stackrel{a}\longrightarrow P'}{\rec(x)P(x)
  \stackrel{a}\longrightarrow P'} &\\
\end{tabular}}
\end{center}
\caption{Operational semantics of process algebra with synchronization algebra $\sigma$}
\label{opsem}
\end{table}

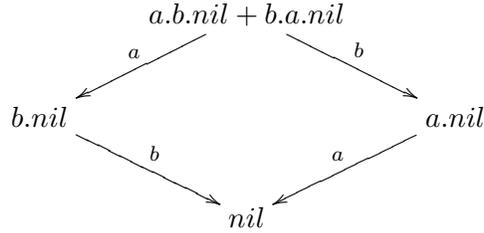
\begin{figure}
\[
\xymatrix{ 
& a.b.nil+b.a.nil \ar@{->}[dl]_-{a} \ar@{->}[dr]^-{b}& \\ b.nil
\ar@{->}[dr]^-{b}&& a.nil \ar@{->}[dl]_-{a}\\ & nil & }
\] 
\caption{Labelled transition system of $a.b.nil+b.a.nil$}
\label{seq}
\end{figure}

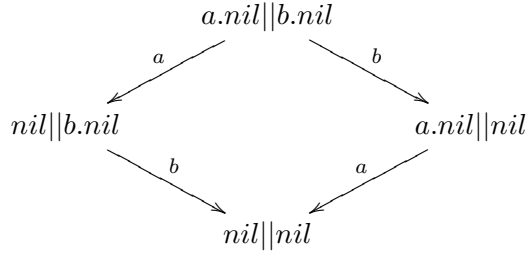
\begin{figure}
\[
\xymatrix{ 
& a.nil|| b.nil \ar@{->}[dl]_-{a} \ar@{->}[dr]^-{b}& \\
nil|| b.nil \ar@{->}[dr]^-{b}&& a.nil|| nil \ar@{->}[dl]_-{a}\\
& nil|| nil & 
}
\]
\caption{Labelled transition system of $a.nil|| b.nil$ with $\sigma(a,b)=\bot$}
\label{conc}
\end{figure}

The operational rules allow us to construct the labelled transition
system dec\-o\-rated by process names associated with a given process
name. For example, consider $P = a.b.nil + b.a.nil \in \proc_\Sigma$.
The Act rule of Table~\ref{opsem} tells us that there exists a
transition $a.b.nil \stackrel{a} \longrightarrow b.nil$ (apply the Act
rule to $\mu=a$ and $P=b.nil$).  Therefore by the Sum1 rule, there
exists a transition $a.b.nil + b.a.nil \stackrel{a}\longrightarrow
b.nil$. The Act rule provides the transition $b.nil\stackrel{b}
\longrightarrow nil$. Symmetrically, one obtains the two other
transitions of Figure~\ref{seq}.  Figure~\ref{conc} describes the
labelled transition system decorated by process names associated with
$a.nil|| b.nil$. The decoration is different from that of $a.b.nil +
b.a.nil$, but the $1$-dimensional paths are the same.

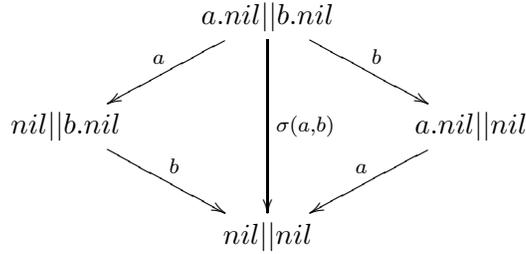
\begin{figure}
\[
\xymatrix{
  & a.nil|| b.nil \ar@{->}[dl]_-{a} \ar@{->}[dd]^-{\sigma(a,b)} \ar@{->}[dr]^-{b}& \\
  nil|| b.nil \ar@{->}[dr]^-{b}& & a.nil|| nil \ar@{->}[dl]_-{a}\\
  & nil|| nil & }
\]
\caption{Labelled transition system of $a.nil|| b.nil$ with synchronization of $a$ and $b$}
\label{com}
\end{figure}

Let $P(x) = \mu.x$. Then the labelled decorated transition system
associated with $\rec(x)P(x)$ is
$\rec(x)P(x)\stackrel{\mu}\longrightarrow \rec(x)P(x)$. Indeed, the
Act rule provides the transition $P(\rec(x)P(x)) \stackrel{\mu}
\longrightarrow \rec(x)P(x)$. Then the Rec rule gives the transition
$\rec(x)P(x) \stackrel{\mu} \longrightarrow
\rec(x)P(x)$. Figure~\ref{com} gives an example of synchronization
obtained by using the Com rule.  Note that in all these examples,
there is a unique initial state which is canonically decorated by the
name of the process we are studying.

Everything is standard in the presentation above except the choice to
consider only acyclic transition systems, which entails the labelling
of the states with processes as in the synchronization tree
semantics. From a directed algebraic topological point of view, this
is equivalent to saying that we consider only universal directed
coverings (dicoverings) as defined in~\cite{MR1994938} or
in~\cite{FR}.

\section{Decorated $\sigma$-labelled precubical set}
\label{slabel}

\subsection{Definition}

\begin{nota}
  Let $[0] = \{()\}$ and $[n] = \{0,1\}^n$ for $n \geq 1$.  By
  convention, one has $\{0,1\}^0=[0]=\{()\}$. The set $[n]$ is
  equipped with the product ordering $\{0<1\}^n$.
\end{nota}

Let $\delta_i^\alpha \colon [n-1] \rightarrow [n]$ be the set map
defined for $1\leq i\leq n$ and $\alpha \in \{0,1\}$ by
\[\delta_i^\alpha(\epsilon_1, \dots, \epsilon_{n-1}) = (\epsilon_1,
\dots, \epsilon_{i-1}, \alpha, \epsilon_i, \dots, \epsilon_{n-1}).\]
The small category $\square$ is by definition the subcategory of the
category of sets with the set of objects $\{[n],n\geq 0\}$ and
generated by the morphisms $\delta_i^\alpha$.  They satisfy the
cocubical relations $\delta_j^\beta \delta_i^\alpha = \delta_i^\alpha
\delta_{j-1}^\beta $ for $i<j$ and for all $(\alpha,\beta)\in
\{0,1\}^2$. If $p>q\geq 0$, then the set of morphisms
$\square([p],[q])$ is empty. If $p = q$, then the set
$\square([p],[p])$ is the singleton $\{\id_{[p]}\}$. For $0\leq p \leq
q$, all maps of $\square$ from $[p]$ to $[q]$ are one-to-one. The
converse is false. The inclusion from $[1]$ to $[2]$ defined by
$(0)\mapsto (0,0)$ and $(1)\mapsto (1,1)$ is not a morphism of
$\square$. Indeed, the category $\square$ does not contain the
compositions of cubes.

A good reference for presheaves is~\cite{MR1300636}.

\bd[\cite{Brown_cube}] The category of presheaves over $\square$,
denoted by $\square^{op}\set$, is called the category of
\emph{precubical sets}.  A precubical set $K$ consists in a family of
sets $(K_n)_{n \geq 0}$ and of set maps $\de_i^\alpha\colon K_n
\rightarrow K_{n-1}$ with $1\leq i \leq n$ and $\alpha\in\{0,1\}$
satisfying the cubical relations $\de_i^\alpha\de_j^\beta =
\de_{j-1}^\beta \de_i^\alpha$ for any $\alpha,\beta\in \{0,1\}$ and
for $i<j$. An element of $K_n$ is called an \emph{$n$-cube}. \ed

Let $\square[n]:=\square(-,[n])$. By Yoneda's lemma, one has the
natural bijection of sets \[\square^{op}\set(\square[n],K)\iso K_n\]
for every precubical set $K$. The boundary of $\square[n]$ is the
precubical set denoted by $\de \square[n]$ defined by removing the
interior of $\square[n]$:
\begin{itemize} 
\item $(\de \square[n])_k := (\square[n])_k$ for $k<n$,
\item $(\de \square[n])_k = \varnothing$ for $k\geq n$.
\end{itemize} 
In particular, one has $\de \square[0] = \varnothing$.

\begin{nota} Let $K$ be a precubical set. Let $K_{\leq n}$ be the
  precubical set obtained from $K$ by keeping the $p$-dimensional
  cubes of $K$ only for $p\leq n$. In particular, $K_{\leq 0}=K_0$.
\end{nota}

\begin{nota} Let $f\colon K \rightarrow L$ be a morphism of precubical
  sets.  Let $n\geq 0$.  The set map from $K_n$ to $L_n$ induced by
  $f$ will be sometimes denoted by $f_n$. \end{nota}

\bd Let $\square_n\subset \square$ be the full subcategory of
$\square$ whose set of objects is $\{[k],k\leq n\}$. The category of
presheaves over $\square_n$ is denoted by $\square_n^{op}\set$. Its
objects are called the \emph{$n$-dimensional precubical sets}.  \ed

We recall now Goubault's construction (in fact a variant of Goubault's
construction) of the precubical set of labels for a trivial
synchronization algebra:

\bp[\cite{labelled}\label{cubeetiquette}] Let
\begin{itemize}
\item $(!\Sigma)_0=\{()\}$ (the empty word),
\item for $n\geq 1$, $(!\Sigma)_n=\Sigma^n$,
\item $\de_i^0(a_1,\dots,a_n) = \de_i^1(a_1,\dots,a_n) =
  (a_1,\dots,\widehat{a_i},\dots,a_n)$, where $\widehat{a_i}$ means
  that $a_i$ is removed.
\end{itemize}
Then these data generate a precubical set $!\Sigma$. \ep

\begin{nota} \label{cubeetiquette1} Let $!^\sigma\Sigma$ be the
  precubical subset of $!\Sigma$ containing the $n$-uples
  $(a_1,\dots,a_n)$ such that $\sigma(a_i,0)=a_i$ for all
  $i$. \end{nota}

So an $n$-cube $(a_1,\dots,a_n)$ of $!\Sigma$ belongs to
$!^\sigma\Sigma$ if and only if the corresponding actions
$a_1,\dots,a_n$ may run asynchronously for the synchronization algebra
$\sigma$.

\begin{figure}
\[
\xymatrix{
& () \ar@{->}[rd]^{(b)}&\\
()\ar@{->}[ru]^{(a)}\ar@{->}[rd]_{(b)} & (a,b) & ()\\
&()\ar@{->}[ru]_{(a)}&}
\]
\caption{Concurrent execution of $a$ and $b$}
\label{concab}
\end{figure}
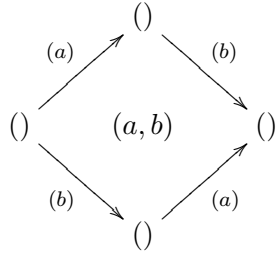

\bd A \emph{$\sigma$-labelled precubical set} is an object of the
comma category \[\square^{op}\set \downarrow !^\sigma\Sigma.\] That
is, an object is a map of precubical sets $\ell\colon K \rightarrow
!^\sigma\Sigma$ and a morphism is a commutative diagram \[ \xymatrix{
  K \ar@{->}[rr]\ar@{->}[rd]&& L \ar@{->}[ld]\\ & !^\sigma\Sigma.&}
\]
The $\ell$ map is called the \emph{labelling map}.  The precubical
set $K$ is sometimes called the \emph{underlying precubical set} of
the $\sigma$-labelled precubical set. \ed

\begin{nota} A $\sigma$-labelled precubical set $K \rightarrow
  !^\sigma\Sigma$ is sometimes denoted by $K$ without explicitly
  mentioning the labelling map.
\end{nota}

Figure~\ref{concab} gives an example of a $\sigma$-labelled
$2$-cube. It represents the concurrent execution of $a$ and $b$. It is
important to notice that two opposite faces of Figure~\ref{concab}
have the same label.

\bd A $\sigma$-labelled precubical set $\ell\colon K\rightarrow
!^\sigma\Sigma$ \emph{decorated by process names} is a
$\sigma$-labelled precubical set together with a set map $d\colon K_0
\rightarrow \proc_\Sigma$ called the \emph{decoration}. \ed

The category of $\sigma$-labelled precubical sets is complete and
cocomplete. This ensures the existence of pullbacks, pushouts and
binary coproducts. All these constructions will be used in the
definition of the denotational semantics of process algebra using
$\sigma$-labelled precubical sets.

We conclude this presentation with useful notation for the sequel:

\begin{nota} Let $K$ and $L$ be two $1$-dimensional $\sigma$-labelled
  precubical sets. Let us denote by $K\p_\sigma L$ the following
  $1$-dimensional $\sigma$-labelled precubical set:
\begin{itemize}
\item $(K\p_\sigma L)_0 = K_0\p L_0$,
\item $(K\p_\sigma L)_1 = (K_1\p L_0) \sqcup (K_0\p L_1) \sqcup
  \{(x,y)\in K_1\p L_1, \sigma(x,y)\neq \bot\}$,
\item $\de_1^\alpha(x,y) = (\de_1^\alpha(x),y)$ for any $(x,y)\in
  K_1\p L_0$,
\item $\de_1^\alpha(x,y) = (x,\de_1^\alpha(y))$ for any $(x,y)\in
  K_0\p L_1$,
\item $\de_1^\alpha(x,y) = (\de_1^\alpha(x),\de_1^\alpha(y))$ for any $(x,y)\in
  K_1\p L_1$,
\item $\ell(x,y)=\ell(x)$ for any $(x,y)\in K_1\p L_0$,
\item $\ell(x,y)=\ell(y)$ for any $(x,y)\in K_0\p L_1$,
\item $\ell(x,y)=\sigma(x,y)$ for any $(x,y)\in K_1\p L_1$ with
  $\sigma(x,y)\neq \bot$.
\end{itemize}
\end{nota}

If $K$ and $L$ are two labelled transition systems associated with two
processes $P$ and $Q$, then the labelled transition system $K\p_\sigma
L$ is the one associated with the process $P|| Q$ by the operational
rules of Table~\ref{opsem}. Since the synchronization algebra $\sigma$
is commutative and associative, one has the natural isomorphisms
$K\p_\sigma L\iso L\p_\sigma K$ and $K\p_\sigma(L\p_\sigma M) \iso
(K\p_\sigma L)\p_\sigma M$ for all $1$-dimensional $\sigma$-labelled
precubical sets $K$, $L$ and $M$.

\subsection{$\sigma$-labelled coskeleton functor}

We present in the framework of precubical sets the construction
presented by Worytkiewicz in the framework of cubical 
sets~\cite{exHDA}, with some slight modifications in the presentation. See
also~\cite[Proposition 3.11]{PhD_Ulrich}. It is discussed in detail for
motivating our $\sigma$-labelled directed coskeleton construction.

\bd \label{def_shell} Let $\ell\colon K\rightarrow !^\sigma\Sigma$ be
a $\sigma$-labelled precubical set. Let $n\geq 1$. A 
\emph{$\sigma$-labelled $n$-shell (or $n$-dimensional shell)} of $K$
is a commutative diagram of precubical sets
\[
\xymatrix{
\de\square[n+1] \ar@{->}[rr]^-{x} \fd{} && K \fd{\ell} \\
&& \\
\square[n+1] \fr{} && !^\sigma\Sigma.}
\]
\ed

Note that if a map $\de\square[n+1] \rightarrow !^\sigma\Sigma$ with
$n\geq 1$ factors as a composite $\de\square[n+1] \rightarrow
\square[n+1] \rightarrow !^\sigma\Sigma$, then this factorization is
unique.

\bp Let $n\geq 1$. The restriction functor $K\mapsto K_{\leq n}$ from
$\square^{op}\set\downarrow !^\sigma\Sigma$ to
$\square_n^{op}\set\downarrow !^\sigma\Sigma$ has a right
adjoint. \ep

\bpf Let us repeat that the proof is easy and that it is given only
for motivating the $\sigma$-labelled directed coskeleton construction
and for introducing some notation.

First of all, let us consider the functor from $\square^{op}_{n+1}\set
\downarrow !^\sigma\Sigma$ to $\square^{op}_{n}\set \downarrow
!^\sigma\Sigma$ induced by the mapping $K\mapsto K_{\leq n}$. Let
$\cosk^\sigma_{n,n+1} \colon \square^{op}_{n}\set \downarrow
!^\sigma\Sigma \rightarrow \square^{op}_{n+1}\set \downarrow
!^\sigma\Sigma$ be the functor defined by:
\[
\xymatrix{
\bigsqcup\limits_{\hbox{labelled $n$-shells}}\de\square[n+1]  \fr{} \fd{}&& K \fd{}\\
&& \\
\bigsqcup\limits_{\hbox{labelled $n$-shells}}\square[n+1]\fr{}&& \cosk^\sigma_{n,n+1}(K). \cocartesien 
}
\] 
Let $\phi\colon K_{\leq n} \rightarrow L$ be a map of $n$-dimensional
$\sigma$-labelled precubical sets. For any $x\in K_{n+1}$, let $\de
x\colon \de\square[n+1] \subset \square[n+1] \rightarrow K$ be the
corresponding morphism of $\sigma$-labelled precubical sets by the
Yoneda lemma. Then consider the commutative diagram of
$\sigma$-labelled precubical sets
\[
\xymatrix{ \de\square[n+1] \fr{\de x} \fd{} && K_{\leq n}
  \fr{\phi}\fd{} && L \fd{}
  \\
  && &&\\
  \square[n+1] \fr{} && !^\sigma\Sigma
  \fr{} && !^\sigma\Sigma.}
\] 
This family of diagrams for $x$ running over $K_{n+1}$ gives rise to a
map of $(n+1)$-dimensional $\sigma$-labelled precubical sets from $K$
to $\cosk^\sigma_{n,n+1}(L)$. Hence the
bijection \[\square^{op}_{n}\set(K_{\leq n},L) \iso
\square^{op}_{n+1}\set(K,\cosk^\sigma_{n,n+1}(L)).\]

Now take a general $\sigma$-labelled precubical set $K\in
\square^{op}\set\downarrow !^\sigma\Sigma$. Let
\[
\cosk^\sigma_{n,n}=\id
\]
and
\[
\cosk^\sigma_{n,n+p}=\cosk^\sigma_{n+p-1,n+p}\circ \dots \circ
\cosk^\sigma_{n,n+1}.
\]
Then the preceding construction gives the
commutative diagram of $\sigma$-labelled precubical sets
\[
\xymatrix{
K_{\leq n} \fr{}\fd{} && K_{\leq n+1} \fr{}\fd{} && K_{\leq n+2} \fr{}\fd{} && \dots \\
&&&&&&\\
\cosk^\sigma_{n,n}(L) \fr{}&&\cosk^\sigma_{n,n+1}(L)
\fr{}&&\cosk^\sigma_{n,n+2}(L) \fr{} &&\dots\, .}
\] 
Hence a map of $\sigma$-labelled precubical sets \[K = \liminj_{k\geq
  0} K_{\leq n+k} \rightarrow \liminj_{k\geq 0}
\cosk^\sigma_{n,n+k}(L).\] So the functor $\liminj_{k\geq 0}
\cosk^\sigma_{n,n+k}$ is the right adjoint. \epf

\begin{nota} Let $\cosk^\sigma_n := \liminj_{k\geq 0}
  \cosk^\sigma_{n,n+k}$. \end{nota}

\bp \label{pblocal} Let $p\geq 2$. Let $\square[p]$ be a
$\sigma$-labelled $p$-cube. Then the $\sigma$-labelled $p$-cube
$\square[p]$ is strictly included in the $\sigma$-labelled precubical
set $\cosk^\sigma_1(\square[p]_{\leq 1})$. \ep

\bpf Let $f(\epsilon_1,\epsilon_2) = (\epsilon_2,\epsilon_1,0,\dots
,0)$ be a set map from $[2]$ to $[p]$. One can consider the
commutative diagram of $\sigma$-labelled precubical sets
\[
\xymatrix{
\de\square[2] \ar@{->}[rr]^-{f} \fd{} && \square[p]_{\leq 1} \fd{\ell} \\
&& \\
\square[2] \fr{} && !^\sigma\Sigma.}
\]
This defines a $2$-cube of $\cosk_1^\sigma(\square[p]_{\leq 1})$. This
$2$-cube does not belong to $\square[p]_2$ since $f$ is not a morphism
of the small category $\square$.  \epf

\begin{figure}
\[
\xymatrix{
&  ()\ar@{->}[rd]^{(b)}&\\
()\ar@{->}[ru]^{(a)}\ar@{->}[rd]_{(b)} & & ()\\
&()\ar@{->}[ru]_{(a)}&}
\]
\caption{Sequential execution of $a$ and $b$ or of $b$ and $a$}
\label{seqab0}
\end{figure}
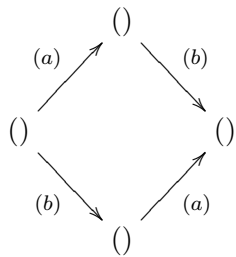

For example, the $1$-dimensional $\sigma$-labelled precubical set
$\square[2]_{\leq 1}$ of Figure~\ref{seqab0} has (with $\sigma(a,0)=a$
and $\sigma(b,0)=b$):
\begin{itemize}
\item two non-degenerate $\sigma$-labelled $1$-dimensional shells
  corresponding to the two set maps from $[2]$ to itself defined by
  $(\epsilon_1,\epsilon_2) \mapsto (\epsilon_1,\epsilon_2)$ and
  $(\epsilon_1,\epsilon_2) \mapsto (\epsilon_2,\epsilon_1)$.
\item If $a = b$, then two degenerate $\sigma$-labelled $1$-dimensional
  shells corresponding to the two set maps $f,g\colon [2]\rightarrow [2]$
  defined by $f(\epsilon_1,\epsilon_2) =
  (\min(\epsilon_1,\epsilon_2),\max(\epsilon_1,\epsilon_2))$ and
  $f(\epsilon_1,\epsilon_2) =
  (\max(\epsilon_1,\epsilon_2),\min(\epsilon_1,\epsilon_2))$. The
  condition $a=b$ comes from the fact that two opposite faces must be
  labelled in the same way.
\end{itemize}

Proposition~\ref{pblocal} means that in the $\sigma$-labelled
coskeleton several different cubes may correspond to the same set of
concurrent actions.

\subsection{$\sigma$-labelled directed coskeleton construction}

We now present the \emph{$\sigma$-labelled directed coskeleton
  construction}. It is defined for any $\sigma$-labelled
$1$-dimensional precubical set $K$ such that $K_0=[p]$ for some $p\geq
2$, as a $\sigma$-labelled precubical subset of $\cosk_1^\sigma(K)$.

The role of the $\sigma$-labelled directed coskeleton is to construct,
from any $\sigma$-labelled precubical set of the form $K\p_\sigma L$
with $K=\square[p]_{\leq 1}$ and $L=\square[q]_{\leq 1}$ for some
$p,q\geq 0$, another $\sigma$-labelled precubical set with the same
$0$-cubes and the same $1$-cubes so that each set of $n$ labelled
$1$-transitions running concurrently is assembled to an $n$-cube
\emph{in exactly one way}.

First of all, we give three examples for motivating the technical
definition.

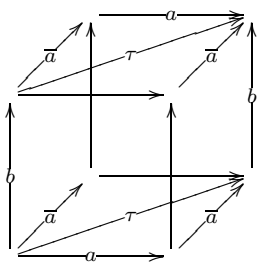
\begin{figure}
\[
\xymatrix{
&\ar@{->}[rr]|{a} &&\\
\ar@{->}[rrru]|-{\tau}\ar@{->}[ru]|-{\overline{a}}\ar@{->}[rr]&&\ar@{->}[ru]|-{\overline{a}}&\\
&\ar@{->}[uu] \ar@{->}[rr]&&\ar@{->}[uu]|-{b}\\
\ar@{->}[uu]|-{b}\ar@{->}[ru]|-{\overline{a}}\ar@{->}[rrru]|-{\tau}\ar@{->}[rr]|-{a}&&\ar@{->}[uu]\ar@{->}[ru]|-{\overline{a}}&
}\]
\caption{1-dimensional paths of $(a.nil|| b.nil)|| \overline{a}.nil$}
\label{3sync00}
\end{figure}

Consider the CCS process $P = (a.nil|| b.nil)|| \overline{a}.nil$ of
Figure~\ref{3sync00}. It corresponds to the concurrent execution with
possible synchronization of a full labelled $2$-cube with a labelled
$1$-cube. In this situation, the actions $a$ and $\overline{a}$ may
synchronize to give an action relabelled by $\tau$, and by definition
of CCS, the actions $\tau$ and $b$ may run concurrently. The
construction of the interpretation $\square|[P|]$ of $P$ starts from
the $1$-dimensional labelled precubical set given by the
$1$-dimensional operational semantics of Table~\ref{opsem}. We then
have to fill all labelled shells corresponding to the six maps of
$\square([2],[3])$ and to the unique map of $\square([3],[3])$.
Moreover, we have to fill the labelled $1$-dimensional shell
corresponding to the set map $f\colon [2] \rightarrow [3]$ defined by
$f(\epsilon_1,\epsilon_2) = (\epsilon_1,\epsilon_2,\epsilon_1)$
corresponding to the concurrent execution of $\tau$ (synchronized
action of $a$ and $\overline{a}$) and $b$. Three important remarks
must be made:
\begin{itemize}
\item The map $f$ is not a map of the small category $\square$.
\item The map $g(\epsilon_1,\epsilon_2) =
  (\epsilon_2,\epsilon_1,\epsilon_2)$ corresponds to the concurrent
  execution of $\tau$ (synchronized action of $a$ and $\overline{a}$)
  and $b$ as well, yet $f\neq g$.
\item The map $g$ can be ruled out using the following fact: the first
  appearance of $\epsilon_2$ is before the first appearance of
  $\epsilon_1$ by reading $(\epsilon_2,\epsilon_1,\epsilon_2)$ from
  the left to the right. The map $g$ is said to be \emph{twisted}. The
  map $f$ is said to be \emph{non-twisted} since the first appearance
  of $\epsilon_1$ is before the first appearance of $\epsilon_2$.
\end{itemize}
The twisted map $g$ corresponds to a $2$-cube of
$\cosk_1^\sigma(\square|[P|]_{\leq 1})$ which does not belong to
$\square|[P|]$.  So there is a strict inclusion of labelled precubical
sets \[\square|[P|] \subset \cosk_1^\sigma(\square|[P|]_{\leq 1}).\]

Consider the CCS process $Q =
(a.nil|| b.nil)|| (\overline{b}.nil|| \overline{a}.nil|| c.nil)$. It
corresponds to the concurrent execution with possible synchronizations
of a full labelled $2$-cube and a full labelled $3$-cube. The labelled
precubical set $\square|[Q|]$ interpreting $Q$ is constructed as
follows:
\begin{itemize}
\item Start from the labelled precubical set $\square[5]_{\leq 1}$ and
  add all $1$-dimensional labelled cubes corresponding to the possible
  synchronizations of $a$ and $\overline{a}$, and of $b$ and
  $\overline{b}$: see the $1$-dimensional operational semantics of
  Table~\ref{opsem}.
\item Add a $(n+1)$-cube $\square[n+1]$ for each map of
  $\square([n+1],[5])$ with $n\geq 1$.
\item We must treat the possible synchronizations of $a$ and
  $\overline{a}$ and of $b$ and $\overline{b}$, each one giving an
  action relabelled by $\tau$. For example, the non-twisted map
  $h_1\colon (\epsilon_1,\epsilon_2)\mapsto
  (\epsilon_1,\epsilon_2,\epsilon_2,0,1)$ corresponds to the
  concurrent execution of $a$ and the action synchronizing $b$ and
  $\overline{b}$, with the action $\overline{a}$ not yet started and
  the action $c$ finished. The non-twisted map $h_1$ will correspond
  to a $2$-cube of $\square|[Q|]$. Another example: the non-twisted
  map $h_2\colon (\epsilon_1,\epsilon_2,\epsilon_3)\mapsto
  (\epsilon_1,\epsilon_2,\epsilon_3,\epsilon_1,1)$ corresponds to the
  concurrent execution of $b$, $\overline{b}$ (which do not
  synchronize here), and the action synchronizing $a$ and
  $\overline{a}$, with the action $c$ finished. The non-twisted map
  $h_2$ will correspond to a $3$-cube of $\square|[Q|]$,
  etc.
\end{itemize} 

Once again, the labelled precubical set $\square|[Q|]$ is strictly
included in the labelled precubical set
$\cosk_1^\sigma(\square|[Q|]_{\leq 1})$.

Consider now the CCS process $R=(a.nil)|| (\overline{a}.nil ||
\overline{a}.nil)$. In this situation, the action $a$ may synchronize
with the left-hand action $a$, or with the right-hand one.  In this
case, the non-twisted mapping $(\epsilon_1) \mapsto
(\epsilon_1,\epsilon_1,0)$ corresponds to the execution of an action
synchronizing $a$ and the left-hand $\overline{a}$, with the
right-hand $\overline{a}$ not yet started. The non-twisted mapping
$(\epsilon_1,\epsilon_2) \mapsto (\epsilon_1,\epsilon_2,\epsilon_1)$
corresponds to the concurrent execution of the left-hand
$\overline{a}$ and the action synchronizing $a$ and the right-hand
$\overline{a}$, etc.

Let us give now the definition of a non-twisted labelled shell.

\bd \label{def_preshell} Let $\ell\colon K\rightarrow !^\sigma\Sigma$
be a $\sigma$-labelled precubical set. Suppose moreover that $K_0=[p]$
for some $p\geq 2$.  Let $n\geq 1$. A \emph{$\sigma$-labelled
  non-twisted $n$-shell (or $n$-dimensional shell)} of $K$ is a
commutative diagram of precubical sets
\[
\xymatrix{
\de\square[n+1] \ar@{->}[rr]^-{x} \fd{} && K \fd{\ell} \\
&& \\
\square[n+1] \fr{} && !^\sigma\Sigma}
\]
such that the set map $x_0\colon [n+1] = \de\square[n+1]_0 \rightarrow [p] =
K_0$ is \emph{non-twisted}.  That is, the set map $x_0\colon [n+1] =
\de\square[n+1]_0 \rightarrow [p] = K_0$ is a composite\footnote{The
  factorization is necessarily unique.}
\[x_0\colon [n+1] \stackrel{\phi}\longrightarrow [q]
\stackrel{\psi}\longrightarrow [p],\] where $\psi$ is a morphism of the
small category $\square$ and where $\phi$ is of the form
\[(\epsilon_1,\dots,\epsilon_{n+1}) \mapsto
(\epsilon_{i_1},\dots,\epsilon_{i_q})\] such that
$\{1,\dots,n+1\}\subset \{i_1,\dots,i_q\}$ and such that the first
appearance of $\epsilon_i$ is before the first appearance of
$\epsilon_{i+1}$ in $ (\epsilon_{i_1},\dots,\epsilon_{i_q})$ for any
$1\leq i\leq n$ by reading from the left to the right.  \ed

The map $\phi$ is not necessarily a morphism of the small category
$\square$. For example, $\phi\colon [3] \rightarrow [5]$ defined by
$\phi(\epsilon_1,\epsilon_2,\epsilon_3)=(\epsilon_1, \epsilon_1,
\epsilon_2, \epsilon_3, \epsilon_3)$ is not a morphism of $\square$.
Note that the set map $x_0$ is then one-to-one.

Let $K$ be an object of $\square_{1}^{op}\set \downarrow
!^\sigma\Sigma$ such that $K_0=[p]$ for some $p\geq 2$.  Let $K^{(n)}$
be the object of $\square_{n}^{op}\set \downarrow !^\sigma\Sigma$
inductively defined for $n\geq 1$ by $K^{(1)} = K$ and by the
following pushout diagram of $\sigma$-labelled precubical sets:
\[
\xymatrix{
\bigsqcup\limits_{\hbox{\emph{non-twisted} $\sigma$-labelled $n$-shells}}\de\square[n+1]  \fr{} \fd{}&& K^{(n)} \fd{}\\
&& \\
\bigsqcup\limits_{\hbox{\emph{non-twisted} $\sigma$-labelled $n$-shells}}\square[n+1]\fr{}&& K^{(n+1)}. \cocartesien 
}
\]
Since $(\de\square[n+1])_p = (\square[n+1])_p$ for $p\leq n$, one has
$(K^{(n+1)})_p = (K^{(n)})_p$ for $p\leq n$, and by construction,
$(K^{(n+1)})_{n+1}$ is the set of non-twisted $\sigma$-labelled
$n$-shells of $K$.  There is an inclusion map $K^{(n)} \rightarrow
K^{(n+1)}$.

\begin{nota} Let $K$ be a $1$-dimensional $\sigma$-labelled precubical
  set with $K_0=[p]$ for some $p\geq 2$.  Then let
    \[\COSK^\sigma(K):= \liminj_{n\geq 1} K^{(n)}.\]
\end{nota}

\bd The $\sigma$-labelled precubical set $\COSK^\sigma(K)$ is called
the \emph{$\sigma$-labelled direct\-ed coskeleton} of $K$. \ed

By construction, the $\sigma$-labelled precubical set
$\COSK^\sigma(K)$ is a $\sigma$-labelled precubical subset of
$\cosk^\sigma_1(K)$.

The construction $\COSK^\sigma$ is not functorial since it is not
defined for all $\sigma$-labelled precubical sets. However one has:

\bp \label{fon} Let $f\colon K \longrightarrow L$ be a morphism of
precubical sets with $K_0=[p_1]$, $L_0=[p_2]$ and such that $f_0\colon [p_1]
\rightarrow [p_2]$ is a morphism of the small category $\square$. Let
$n\geq 1$, and let
\[
\xymatrix{
\de\square[n+1] \ar@{->}[rr]^-{x} \fd{} && K \fd{\ell} \\
&& \\
\square[n+1] \fr{} && !^\sigma\Sigma}
\]
be a $\sigma$-labelled non-twisted shell of $K$. Then the commutative
square
\[
\xymatrix{
\de\square[n+1] \ar@{->}[rr]^-{f\circ x} \fd{} && L \fd{\ell} \\
&& \\
\square[n+1] \fr{} && !^\sigma\Sigma}
\]
is a $\sigma$-labelled non-twisted shell of $L$.  \ep

The hypothesis about $f_0$ is necessary. Indeed, take $K = K_0 =
[p_1]$ and $L = L_0 = [p_2]$. Then a morphism of precubical sets
$f\colon K \rightarrow L$ can be any set map from $[p_1]$ to $[p_2]$.

\bpf
The map $x_0\colon [n+1] \rightarrow [p_1]$ is a composite 
\[x_0\colon [n+1] \stackrel{\phi}\longrightarrow [q]
\stackrel{\psi}\longrightarrow [p_1],\] where $\psi$ is a morphism of
the small category $\square$ and where $\phi$ is of the form
\[(\epsilon_1,\dots,\epsilon_{n+1}) \mapsto
(\epsilon_{i_1},\dots,\epsilon_{i_q})\] such that
$\{1,\dots,n+1\}\subset \{i_1,\dots,i_q\}$ and such that the first
appearance of $\epsilon_i$ is before the first appearance of
$\epsilon_{i+1}$ in $ (\epsilon_{i_1},\dots,\epsilon_{i_q})$ for any
$1\leq i\leq n$ by reading from the left to the right.  So the map
$f_0\circ x_0\colon [n+1] \rightarrow [p_2]$ is the composite
\[f_0\circ x_0\colon [n+1] \stackrel{\phi}\longrightarrow [q]
\stackrel{f_0\circ \psi}\longrightarrow [p_2].\]  \epf

Propositions~\ref{pourquoi_non_degenere} and~\ref{corpourquoi} explain
why this new construction works.

\bp \label{pourquoi_non_degenere} Let $\square[p]$ (resp.
$\square[q]$) be a $\sigma$-labelled cube, which corresponds to the
concurrent execution of $p$ actions $(a_1.nil) || \dots || (a_p.nil)$
(resp.\ of $q$ actions $(a_{p+1}.nil) || \dots || \break (a_{p+q}.nil))$.
Let $(A,B,C^{-},C^{+})$ be a partition of the set 
$\{1,\dots,p+q\}$. Let \[f\colon B\cap\{1,\dots,p\} \rightarrow
B\cap\{p+1,\dots,p+q\}\] be a bijection. Then there exists a unique
non-twisted map $g\colon [r]\rightarrow [p+q]$ corresponding to the
concurrent executions of the actions $a_i$ for $i\in A$ and of the
actions synchronizing $a_i$ and $a_{f(i)}$ for $i\in
B\cap\{1,\dots,p\}$, with the actions $a_i$ for $i\in C^{-}$ not yet
started and the actions of $a_i$ for $i\in C^{+}$ already finished.
\ep

\bpf[Sketch of proof] It may be very helpful for the reader to read
the examples given in the beginning of this section.

One necessarily has $r=|A|+|B\cap\{1,\dots,p\}|$, where $|A|$ (resp.\
$|B\cap\{1,\dots,p\}|$) is the cardinal of $A$ (resp.\
$B\cap\{1,\dots,p\}$). The non-twisted map $g\colon [r]\rightarrow
[p+q]$ is obtained as the composite $[r]\stackrel{\phi}\rightarrow [s]
\stackrel{\psi}\rightarrow [p+q]$ defined as follows.  Let
\[
s=|A|+2|B\cap\{1,\dots,p\}| =
|A|+|B\cap\{1,\dots,p\}|+|B\cap\{p+1,\dots,p+q\}|.
\]
One must have $s+|C^{-}|+|C^{+}| = p+q$. Let $C^{-} \cup C^{+} =
\{a_{i_1}, \dots, a_{i_{p+q-s}}\}$ with $i_1<\dots<i_{p+q-s}$. Then
the morphism of $\psi \colon [s] \rightarrow [p+q]$ of the small
category $\square$ is necessarily equal to
$\delta_{i_{p+q-s}}^{\eta(p+q-s)} \circ \dots
\circ\delta_{i_1}^{\eta(1)}$ with $\eta(k)=0$ if $a_{i_k}\in C^{-}$
and $\eta(k)=1$ if $a_{i_k}\in C^{+}$. The map $\phi\colon
[r]\rightarrow[s]$ is constructed as follows. Write $A\cup
(B\cap\{1,\dots,p\}) = \{j_1<\dots<j_r\}$. Then rewrite the set
$\{j_1,\dots,j_r\} \cup \{f(j_k),j_k\in B\cap\{1,\dots,p\}\}$ in
increasing order and replace each occurrence of $j_k$ and $f(j_k)$ by
$\epsilon_k$. One obtains a word using $\epsilon_1,\dots,\epsilon_r$
with $s-r$ repetitions defining a non-twisted map $\phi\colon
[r]\rightarrow[s]$.  \epf

\bp \label{corpourquoi} Let $\square[p]$ be a $\sigma$-labelled full
$p$-cube with $p\geq 2$. Then one has the isomorphism of
$\sigma$-labelled precubical sets $\COSK^\sigma(\square[p]_{\leq 1})
\iso \square[p]$. \ep

\bpf Let $n\geq 2$, $\square[n] \rightarrow \square[p]$ be a $n$-cube
of $\square[p]$ with $n\geq 2$.  Let $X$ be the composite map $x\colon
\de\square[n] \subset \square[n] \rightarrow \square[p]$.  Thus one
obtains a commutative square
\[
\xymatrix{
\de\square[n] \ar@{->}[rr]^-{x} \fd{} && \square[p] \fd{\ell} \\
&& \\
\square[n] \fr{} && !^\sigma\Sigma}
\]
and therefore a $\sigma$-labelled $(n-1)$-shell. This shell is
non-twisted since the map $x_0$ is the composite $[n]
\stackrel{\id}\rightarrow [n] \stackrel{x_0}\rightarrow [p]$.

Conversely, start from a $\sigma$-labelled non-twisted $(n-1)$-shell
as above with $n\geq 2$. The map $x$ induces a non-twisted map
\[
x_0\colon 
\de\square[n]_0 = \square([0],[n]) = [n] \rightarrow \square([0],[p])
= \square[p]_0 = [p].
\]
The map $x_0$ factors as a composite $[n] \stackrel{\phi}\rightarrow
[m] \stackrel{\psi}\rightarrow [p]$ as in the definition of a
non-twisted shell. The map $\phi$ cannot contain any repetition since
there are no synchronizations by definition of $\square[p]$. So $\phi$
is the identity of $[n]=[m]$ and $x_0=\psi$ is a morphism of the small
category $\square$.

Thus there exists a bijective correspondence between the
$\sigma$-labelled $p$-cubes of $\square[n]$ and the non-twisted
$\sigma$-labelled $(p-1)$-shells of $\square[n]$ for $p\geq 2$.  Hence
we have the result.  \epf

\section{Denotational semantics using $\sigma$-labelled precubical sets}
\label{semccsprecube}

It is defined by induction on the syntax a $\sigma$-labelled
precubical set $\square|[P|]$ for each process name
$P\in\proc_\Sigma$. The $\sigma$-labelled precubical set
$\square|[P|]$ will have a unique initial state canonically decorated
by the process name $P$ and its other states will be decorated as well
in an inductive way. Therefore for every process name $P$,
$\square|[P|]$ will be an object of $\{i\}\downarrow \square^{op}\set
\downarrow !^\sigma\Sigma$ equipped with a decoration by process
names. The only new interpretation is the one of $P\| Q$.  The other
ones are well-known.

\subsection{Interpretation of  $nil$}

The process $nil$ is the idle process. Therefore
\[\boxed{\square|[nil|]:=\square[0]}  \]  

\subsection{Interpretation of $\mu.nil$}

The decorated $\sigma$-labelled precubical set associated with
$\mu.nil$ consists in the unique labelled transition $\mu.nil
\stackrel{\mu} \longrightarrow
nil$. Therefore \[\boxed{\square|[\mu.nil|]:=\mu.nil
  \stackrel{\mu}\longrightarrow nil} \]

\subsection{Interpretation of $\mu.P$}

The Act rule of Table~\ref{opsem} provides the transition $\mu.P
\stackrel{\mu} \longrightarrow P$. The interpretation of $\mu.P$ is
then obtained by identifying the final state of $\square|[\mu.nil|]$
with the initial state of $\square|[P|]$. Thus, by definition, one has
the cocartesian diagram
\[\boxed{\xymatrix{
\square[0]=\{0\} \ar@{->}[rr]^-{0\mapsto nil} \ar@{->}[dd]^-{0\mapsto P} && \square|[\mu.nil|] \ar@{->}[dd] \\
&&\\
\square|[P|] \ar@{->}[rr] && \cocartesien {\square|[\mu.P|]}}}
\]
The mapping $\square|[P|] \mapsto \square|[\mu.P|]$ is functorial with
respect to 
\[
\square|[P|]\in \{i\}\downarrow \square^{op}\set
\downarrow !^\sigma\Sigma.
\]

\subsection{Interpretation of $P+Q$}

The interpretation of $P+Q$ is obtained by a binary coproduct. Thus
\[\boxed{\square|[P+Q|] := \square|[P|] \oplus \square|[Q|]}\] where $\oplus$ is
the binary coproduct in the category $\{i\}\downarrow
\square^{op}\set \downarrow !^\sigma\Sigma$. Indeed, the precubical
set $\square|[P+Q|]$ must have a unique initial state. The
construction of $\square|[P+Q|]$ is functorial with respect to
$\square|[P|]$ and $\square|[Q|]$ in $\{i\}\downarrow
\square^{op}\set \downarrow !^\sigma\Sigma$.

\subsection{Interpretation of $(\nu a)P$}

The restriction rule
\[\inference[\bf{Res}]{P \stackrel{\mu} \longrightarrow P'   \mu\neq a \hbox{ and }\sigma(\mu,a)=\bot}{(\nu a)P \stackrel{\mu} \longrightarrow (\nu a)P'}\] 
tells us that all events which may synchronize with $a$ and also $a$
itself must be removed. Thus, the interpretation of $(\nu a)P$ is
defined by the following pullback diagram:
\[\boxed{
\xymatrix{
  \square|[(\nu a) P|] \ar@{->}[rr] \ar@{->}[dd] \cartesien && \square|[P|] \ar@{->}[dd] \\
  && \\
  !^\sigma(\Sigma\backslash (\{a\}\cup\{b,\sigma(a,b)\neq \bot\})) \ar@{->}[rr]
  && !^\sigma\Sigma}}
\] 
The mapping \[\square|[P|] \mapsto \square|[(\nu a)P|]\] is functorial
with respect to $\square|[P|]\in \{i\}\downarrow \square^{op}\set
\downarrow !^\sigma\Sigma$.

\subsection{Interpretation of $P|| Q$}

The most complicated construction is coming now. One has to construct
the interpretation of $P|| Q$ from the knowledge of the precubical
sets $\square|[P|]$ and $\square|[Q|]$. Let $\square[p] \rightarrow
\square|[P|]$ be a $p$-cube of $\square|[P|]$ corresponding to the
concurrent execution of $p$ transitions $(a_1.nil) || \dots ||
(a_p.nil)$ in $P$. Let $\square[q] \rightarrow \square|[Q|]$ be a
$q$-cube of $\square|[Q|]$ corresponding to the concurrent execution
$(a_{p+1}.nil)|| \dots|| (a_{p+q}.nil)$ of $q$ transitions in $Q$.
First we construct all possible synchronizations between $\square[p]
\rightarrow !^\sigma\Sigma$ and $\square[q] \rightarrow
!^\sigma\Sigma$. The operational semantics of process algebra (cf.\
Figure~\ref{com}) gives us the $1$-dimensional $\sigma$-labelled
precubical set $\square[p]_{\leq 1} \p_\sigma \square[q]_{\leq 1}$.
Then let \[\boxed{\square[p] \ot_\sigma \square[q] :=
  \COSK^\sigma(\square[p]_{\leq 1} \p_\sigma \square[q]_{\leq 1})} \]
The effect of the $\COSK^\sigma$ construction is to add all possible
concurrent executions by filling all $\sigma$-labelled non-twisted
shells of higher dimension. Hence we have the construction of
$\square|[P|| Q|]$ using Proposition~\ref{fon}:
\[\square|[P|| Q|]  := \liminj_{\square[p]\rightarrow \square|[P|]} \liminj_{\square[q]\rightarrow \square|[Q|]}\square[p] \ot_{\sigma} \square[q]\] 
or, more explicitly, 
\[\boxed{\square|[P|| Q|]  := \liminj_{\square[p]\rightarrow
\square|[P|]} \liminj_{\square[q]\rightarrow \square|[Q|]}  \COSK^\sigma
\lp\square[p]_{\leq 1} \p_\sigma \square[q]_{\leq 1}\rp} \]

In particular, one has

\bp For any process name $P$, one has the isomorphisms of
$\sigma$-labelled precubical sets $\square|[P|| |nil|] \iso
\square|[nil|| P|] \iso \square|[P|]$. \ep

Note that the role of Proposition~\ref{corpourquoi} is crucial in the
proof of this proposition.

\bpf One has
\begin{align*}
\square|[P|| nil|] & \iso \liminj_{\square[p]\rightarrow \square|[P|]}
\liminj_{\square[q]\rightarrow \square|[nil|]}  \COSK^\sigma\lp
\square[p]_{\leq 1} \p_\sigma \square[q]_{\leq 1}\rp  \\
&\iso \liminj_{\square[p]\rightarrow \square|[P|]}
  \COSK^\sigma\lp \square[p]_{\leq 1} \p_\sigma \square[0]_{\leq 1}
\rp  \\
&\iso \liminj_{\square[p]\rightarrow \square|[P|]}
  \COSK^\sigma\lp
\square[p]_{\leq 1}\rp  \\
& \iso \liminj_{\square[p]\rightarrow \square|[P|]} \square[p] 
\hbox{ by Proposition~\ref{corpourquoi}}\\
&\iso \square|[P|].\qedhere 
\end{align*}
\epf 

One has the bijection of sets \beas
\square|[P|| Q|]_0 && \iso \liminj_{\square[m]\rightarrow \square|[P|]} \liminj_{\square[n]\rightarrow \square|[Q|]}(\square[m] \ot_{\sigma} \square[n])_0 \\
&&\iso \liminj_{\square[m]\rightarrow \square|[P|]}
\liminj_{\square[n]\rightarrow \square|[Q|]}
(\square[m])_0\p (\square[n])_0 \\
&& \iso \square|[P|]_0 \p \square|[Q|]_0.  \eeas Therefore the
construction of $\square|[P|| Q|]$ is functorial with respect to
$\square|[P|]$ and $\square|[Q|]$ as an object of $\{i\}\downarrow
\square^{op}\set \downarrow !^\sigma\Sigma$.

Hence we have the definition: 

\bd Let $K$ and $L$ be two $\sigma$-labelled precubical sets. The
\emph{tensor product with synchronization} (or \emph{synchronized tensor 
  product}) of $K$ and $L$ is
\[K \ot_\sigma L  := \liminj_{\square[p]\rightarrow K} \liminj_{\square[q]\rightarrow L}\square[p] \ot_{\sigma} \square[q].\]
\ed 

Since the synchronization algebra $\sigma$ is commutative, the
underlying precubical sets of $K \ot_\sigma L$ and $L \ot_\sigma K$
are naturally isomorphic for all $\sigma$-labelled precubical sets $K$
and $L$. One also has the natural isomorphisms $K \ot_\sigma (L
\ot_\sigma M) \iso (K \ot_\sigma L) \ot_\sigma M$ (the proof is
postponed until Appendix~\ref{associativity}), and $K\ot_\sigma
\square[0] \iso \square[0] \ot_\sigma K \iso K$ for any
$\sigma$-labelled precubical set $K$, $L$ and $M$.  The particular
case of the trivial synchronization algebra is interesting to notice:

\bp Let $\sigma=\bot$. Let $K$ and $L$ be two $\sigma$-labelled
precubical sets. Then the tensor product with synchronization
$\ot_\sigma$ is the usual tensor product;  that is, 
\[K \ot_\bot L \iso \liminj_{\square[p] \rightarrow
  K}\liminj_{\square[q] \rightarrow L} \square[p+q].\] \ep

\bpf One has $\square[p]_{\leq 1} \p_{\bot} \square[q]_{\leq 1} \iso
\square[p+q]_{\leq 1}$.  Hence the result by
Proposition~\ref{corpourquoi}. \epf

The restriction in dimension $1$ of the synchronized tensor product is
interesting too:

\bp \label{res1} Let $K$ and $L$ be two $\sigma$-labelled precubical
sets. Then one has 
\[
(K \ot_\sigma L)_{\leq 1} \iso K_{\leq 1}
\p_\sigma  L_{\leq 1}.
\]
\ep 

\bpf The formula $(K \ot_\sigma L)_0=K_0\p L_0$ is proved above. One
has
\[(K \ot_\sigma L)_{1} \iso \liminj_{\square[p] \rightarrow
  K}\liminj_{\square[q] \rightarrow L} (\square[p] \ot_{\sigma}
\square[q])_{1}\] since the functor $K\mapsto K_{1}$ preserves all
colimits. So the set $(K \ot_\sigma L)_{1}$ is equal to
{\small\[\liminj_{\square[p] \rightarrow K}\liminj_{\square[q]
    \rightarrow L} (\square[p]_1 \p \square[q]_0) \oplus
  (\square[p]_0\p \square[q]_1) \oplus \{(x,y)\in \square[p]_1\p
  \square[q]_1,\sigma(x,y)\neq \bot\}.\]} The functors $K\mapsto K_0$
and $K\mapsto K_1$ preserve all colimits and the category of sets is
cartesian closed. So one obtains
\[\liminj_{\square[p] \rightarrow
  K}\liminj_{\square[q] \rightarrow L} (\square[p]_1 \p \square[q]_0)
\iso K_1 \p L_0,\] 
and 
\[\liminj_{\square[p] \rightarrow
  K}\liminj_{\square[q] \rightarrow L} (\square[p]_0 \p \square[q]_1) \iso K_0 \p L_1.\]
Finally, the colimit 
\[\liminj_{\square[p] \rightarrow
  K}\liminj_{\square[q] \rightarrow L} \{(x,y)\in \square[p]_1\p
\square[q]_1,\sigma(x,y)\neq \bot\}\] 
is calculated for each value of $\sigma(x,y)$ and one obtains 
{\small\[\liminj_{\square[p] \rightarrow
  K}\liminj_{\square[q] \rightarrow L} \{(x,y)\in \square[p]_1\p
\square[q]_1,\sigma(x,y)\neq \bot\} \iso \{(x,y)\in K_1\p
L_1,\sigma(x,y)\neq \bot\}.\qedhere\]}
\epf

\subsection{Interpretation of $\rec(x)P(x)$}

The interpretation of $\rec(x)P(x)$ is constructed as the least fixed
point of $P$~\cite{MR1365754}. The functoriality of the previous
constructions implies that the mapping $x\mapsto P(x)$ from
$\proc_\Sigma$ to itself induces a functorial mapping \[\square|[x|]
\mapsto \square|[P(x)|]\] from $\{i\}\downarrow \square^{op}\set
\downarrow !^\sigma\Sigma$ to itself.  This fact can be easily checked
by induction on the syntax of $P(x)$. The $\sigma$-labelled precubical
set $\square|[P(nil)|]$ has a unique initial state canonically
decorated by the process name $P(nil)$. Thus we have a map
$\square|[nil|]=\square[0] \rightarrow \square|[P(nil)|]$, and by
functoriality, a map $\square|[P^n(nil)|] \rightarrow
\square|[P^{n+1}(nil)|]$ of $\{i\}\downarrow \square^{op}\set
\downarrow !^\sigma\Sigma$.  As usual, let
\[\boxed{\square|[\rec(x)P(x)|]:=\liminj_n \square|[P^n(nil)|]} \]

As an example, consider the case $P(x)=\mu.x$. Then
$\square|[P^n(nil)|]$ is the decorated $\sigma$-labelled precubical
set
\[\rec(x)P(x) \stackrel{\mu} \longrightarrow \rec(x)P(x) 
\stackrel{\mu} \longrightarrow \dots \rec(x)P(x) \hbox{ ($n$ times
  $\mu$)}.\] Therefore the decorated $\sigma$-labelled precubical set
associated with $\rec(x)P(x)$ is
\[\rec(x)P(x) \stackrel{\mu} \longrightarrow \rec(x)P(x) 
\stackrel{\mu} \longrightarrow \rec(x)P(x) \stackrel{\mu} \longrightarrow \dots 
\hbox{ (indefinitely)}.\] 

\section{Restriction in dimension $1$ and HDA paradigm}
\label{c0}

\bth \label{restrict1} Let $P$ be a process name. Then the
$\sigma$-labelled $1$-dimensional precubical set $\square|[P|]_{\leq
  1}$ coincides with the labelled transition system given by the
operational semantics of Table~\ref{opsem}. \eth

\bpf By induction on the syntax of $P$ using Proposition~\ref{res1}
and the fact that the functor $K\mapsto K_{\leq 1}$ from precubical
sets to $1$-dimensional precubical sets preserves all limits and all
colimits.  \epf

Theorem~\ref{restrict1} is not a soundness result with respect to any
kind of bisimulation. For example, for any process name $P$ of CCS,
the processes $P+P$ and $P$ are strongly bisimilar, whereas the
corresponding labelled flows $|\square|[P+P|]|$ and $|\square|[P|]|$,
where $| -|$ is the geometric realization functor of Section~\ref{rea},
are not weakly equivalent in the model category of flows. The link
with bisimilarity will be the subject of future papers.

By Proposition~\ref{pourquoi_non_degenere}, this semantics satisfies
the paradigm of higher dimensional automata. More precisely, one has

\bth \label{paradigm} Let $P$ be a process name, and $p\geq
1$. Consider a $\sigma$-labelled $p$-shell of $\square|[P|]$:
\[
\xymatrix{
\de\square[p+1] \ar@{->}[rr]^-{x} \fd{} && \square|[P|] \fd{\ell} \\
&& \\
\square[p+1] \ar@{->}[rruu]^-{k}\fr{} && !^\sigma\Sigma.}
\]
Then there exists at most one lift $k$.  \eth

The existence of the lift $k$ means that the $p+1$ actions
$a_1,\dots,a_{p+1}$ we are considering run concurrently in the process
$P$. The unique lift $k$ then corresponds to a unique $(p+1)$-cube.
This is precisely the paradigm of higher dimensional automata.

\bpf By induction on the syntax of $P$. It is clear that if $P$ and
$Q$ satisfy the statement of the theorem, then $P+Q$ satisfies the
same statement. If $P$ satisfies the statement of the theorem, then
$(\nu a)P$ satisfies the same statement since the precubical set
$\square|[(\nu a)P|]$ is a precubical subset of $\square|[P|]$.  If
$P(y)$ is a process name with one free variable $y$, then the
precubical set $\square|[P^n(nil)|]$ is a precubical subset of
$\square|[P^{n+1}(nil)|]$. So $\rec(y) P(y)$ satisfies the statement
of the theorem. Finally, consider two lifts $k_1,k_2\colon \square[p+1]
\rightarrow \square|[P|| Q|]$. Since one has the colimit of sets
\[\square|[P|| Q|]_{p+1}  := \liminj_{\square[m]\rightarrow \square|[P|]} \liminj_{\square[n]\rightarrow \square|[Q|]}(\square[m] \ot_{\sigma} \square[n])_{p+1},\] 
the maps $k_1$ and $k_2$ factor as composites \[k_1 \colon
\square[p+1] \rightarrow \square[m_1] \ot_{\sigma} \square[n_1]
\rightarrow \square|[P|| Q|],\] and \[k_2 \colon \square[p+1]
\rightarrow \square[m_2] \ot_{\sigma} \square[n_2] \rightarrow
\square|[P|| Q|].\] The subcategory of $(\square\downarrow
\square|[P|]) \p (\square\downarrow \square|[Q|])$ of objects $(m,n)$
such that $x$ factors as a composite $x\colon \de\square[p+1]
\rightarrow \square[m] \ot_{\sigma} \square[n] \rightarrow
\square|[P|| Q|]$ is filtered.  Therefore one can suppose that
$(m_1,n_1) = (m_2,n_2)$. Thus $k_1 = k_2$ by
Proposition~\ref{pourquoi_non_degenere}.  So $P||Q$ satisfies the
statement of the theorem. \epf

Theorem~\ref{paradigm} implies that not all $\sigma$-labelled
precubical sets can be viewed as higher dimensional automata:

\begin{corollary} \label{sfa} Assume that there exist two actions $a$
  and $b$ with \[\sigma(a,0) = a, \sigma(b,0) = b.\]
  Consider the $\sigma$-labelled precubical set $\square[2]_{\leq 1} =
  \de\square[2]$:
\[
\xymatrix{
& () \ar@{->}[rd]^{(b)}&\\ ()\ar@{->}[ru]^{(a)}\ar@{->}[rd]_{(b)} & &
().\\ &()\ar@{->}[ru]_{(a)}&}
\]
Add two $2$-cubes labelled by $(a,b)$. Then the $\sigma$-labelled
precubical set we obtain is not isomorphic to any $\sigma$-labelled
precubical set of the form $\square|[P|]$ for any process name $P$.
\end{corollary}

Corollary~\ref{sfa} is not surprising. Algebraic-topological models of
concurrency are all considerable generalizations of the ``usual''
models of concurrency. This high level of generality is necessary to
obtain convenient settings for doing homotopy.

\section{Decorated $\sigma$-labelled flow}
\label{sflow}

The category $\top$ of \emph{compactly generated topological spaces}
(i.e.\ of weak Hausdorff $k$-spaces) is complete, cocomplete and
cartesian closed (more details for these kinds of topological spaces
are in~\cite{MR90k:54001, MR2000h:55002}, the appendix
of~\cite{Ref_wH} and also in the preliminaries of \cite{model3}). For
the sequel, all topological spaces will be supposed to be compactly
generated. A \emph{compact space} is always Hausdorff.

The category $\top$ is equipped with the unique model structure having
the \emph{weak homotopy equivalences} as weak equivalences and
having the \emph{Serre fibrations}\footnote{that is, a continuous
  map having the RLP with respect to the inclusion $\mathbf{D}^n\p
  0\subset \mathbf{D}^n\p [0,1]$ for any $n\geq 0$ where
  $\mathbf{D}^n$ is the $n$-dimensional disk.} as fibrations.

\bd[\cite{model3}] A \emph{(time) flow} $X$ is a small topological
category without identity maps. The set of objects is denoted by
$X^0$.  The topological space of morphisms from $\alpha$ to $\beta$ is
denoted by $\P_{\alpha,\beta}X$. The elements of $X^0$ are also called
the \emph{states} of $X$. The elements of $\P_{\alpha,\beta}X$ are
called the \emph{(non-constant) execution paths from $\alpha$ to
  $\beta$}. A flow $X$ is \emph{loopless} if for every $\alpha\in
X^0$, the space $\P_{\alpha,\alpha}X$ is empty. \ed

\begin{nota} Let $\P X = \bigsqcup_{(\alpha,\beta)\in X^0\p X^0}
  \P_{\alpha,\beta}X$.  $\P X$ is called the \emph{path space} of $X$.
  The source map (resp.\ the target map) $\P X\rightarrow X^0$ is
  denoted by $s$ (resp.\ $t$).
\end{nota}

\bd Let $X$ be a flow, and let $\alpha \in X^0$ be a state of $X$. The
state $\alpha$ is \emph{initial} if $\alpha\notin t(\P X)$, and the state
$\alpha$ is \emph{final} if $\alpha\notin s(\P X)$. \ed

\bd A morphism of flows $f\colon  X \rightarrow Y$ consists in a set map
$f^0\colon X^0 \rightarrow Y^0$ and a continuous map $\P f\colon \P X \rightarrow
\P Y$ compatible with the structure. The corresponding category is
denoted by $\dtop$. \ed

The strictly associative composition law 
\[
\left\{ \begin{array}{c} \P_{\alpha,\beta}X \p \P_{\beta,\gamma}X \longrightarrow \P_{\alpha,\gamma}X \\
(x,y) \mapsto x*y \end{array} \right.
\]
models the composition of non-constant execution paths. The
composition law $*$ is extended in a usual way to states, that is to
constant execution paths, by $x*t(x) = x$ and $s(x)*x = x$ for every
non-constant execution path $x$.

\begin{nota} 
  The category of sets is denoted by $\set$. The category of par\-tially
  ordered sets or posets together with the strictly increasing maps
  ($x<y$ implies $f(x)<f(y)$) is denoted by $\poset$.
\end{nota}

\begin{figure}
\begin{center}
\includegraphics[width=7cm]{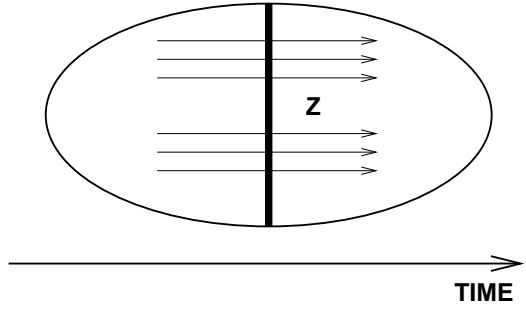}
\end{center}
\caption{Symbolic representation
of $\glob(Z)$ for some topological space $Z$} 
\label{exglob}
\end{figure}

Here are four fundamental examples of flows:
\begin{enumerate}
\item Let $S$ be a set. The flow associated with $S$, still denoted by
  $S$, has $S$ as a set of states and the empty space as path space.
  This construction induces a functor $\set \rightarrow \dtop$ from
  the category of sets to that of flows. The flow associated with a
  set is loopless.
\item Let $(P,\leq)$ be a poset. The flow associated with $(P,\leq)$,
  and still denoted by $P$ is defined as follows: the set of states of
  $P$ is the underlying set of $P$; the space of morphisms from
  $\alpha$ to $\beta$ is empty if $\alpha\geq \beta$ and equals to
  $\{(\alpha,\beta)\}$ if $\alpha<\beta$ and the composition law is
  defined by $(\alpha,\beta)*(\beta,\gamma) = (\alpha,\gamma)$. This
  construction induces a functor $\poset \rightarrow \dtop$ from the
  category of posets together with the strictly increasing maps to the
  category of flows. The flow associated with a poset is loopless.
\item The flow $\glob(Z)$ is defined by
\begin{align*} 
    \glob(Z)^0 &=\{\widehat{0},\widehat{1}\}, \\
    \P \glob(Z) &=  \P_{\widehat{0},\widehat{1}} \glob(Z) = Z, \\
    s &=\widehat{0}, \\
   t &=\widehat{1}
\end{align*}

and a trivial composition law (cf.\
Figure~\ref{exglob}).  It is called the \emph{globe} of $Z$. 
\item The \emph{directed segment} $\vI$ is by definition
  $\glob(\{0\}) \iso \{\widehat{0} < \widehat{1}\}$. 
\end{enumerate}

The model structure of $\dtop$ is characterized as 
follows~\cite{model3}:
\begin{itemize}
\item The weak equivalences are the \emph{weak S-homotopy
    equivalences}, i.e.\ the morphisms of flows $f\colon
  X\longrightarrow Y$ such that $f^0\colon X^0\longrightarrow Y^0$ is
  a bijection of sets and such that $\P f\colon \P X\longrightarrow \P
  Y$ is a weak homotopy equivalence.
\item The fibrations are the morphisms of flows $f\colon
  X\longrightarrow Y$ such that $\P f\colon \P X\longrightarrow \P Y$
  is a Serre fibration.
\end{itemize}
This model structure is cofibrantly generated. The set of generating
cofibrations is the set $I^{gl}_{+}=I^{gl}\cup \{R\colon
\{0,1\}\longrightarrow \{0\},C\colon \varnothing\longrightarrow
\{0\}\}$ with
\[I^{gl}=\{\glob(\mathbf{S}^{n-1})\subset \glob(\mathbf{D}^{n}), n\geq
0\},\] where $\mathbf{D}^{n}$ is the $n$-dimensional disk and
$\mathbf{S}^{n-1}$ the $(n-1)$-dimensional sphere. By convention, the
$(-1)$-dimensional sphere is the empty space. The set of generating
trivial cofibrations is
\[J^{gl}=\{\glob(\mathbf{D}^{n}\p\{0\})\subset
\glob(\mathbf{D}^{n}\p [0,1]), n\geq 0\}.\] 

 \begin{nota} The cofibrant replacement functor is denoted by
  $(-)^{\textit{cof}}$. \end{nota}

\bd The \emph{flow of labels} $?^\sigma\Sigma$ is defined as follows:
$(?^\sigma\Sigma)^0 = \{0\}$ and $\P ?^\sigma\Sigma$ is the discrete
free associative monoid without unit generated by the elements of
$\Sigma$ and by the algebraic relations $a*b=b*a$ if and only if $a$
and $b$ can occur asynchronously; that is, $\sigma(a,0)=a$ and
$\sigma(b,0)=b$.  \ed

\bd A \emph{$\sigma$-labelled flow} is an object of the comma category
$\dtop \downarrow ?^\sigma\Sigma$. That is an object is a map of
flows $\ell\colon X \rightarrow ?^\sigma\Sigma$ and a morphism is a
commutative diagram
\[ \xymatrix{ X \ar@{->}[rr]\ar@{->}[rd]&& Y \ar@{->}[ld]\\ &
  ?^\sigma\Sigma.&}
\]
The $\ell$ map is called the \emph{labelling map}.  The flow $X$ is
sometimes called the \emph{underlying flow} of the $\sigma$-labelled
flow.  \ed

\bd A $\sigma$-labelled flow $\ell\colon X \rightarrow ?^\sigma\Sigma$
\emph{decorated by process names} is a $\sigma$-labelled flow together
with a set map $d\colon X^0 \rightarrow \proc_\Sigma$ called the 
 \emph{decoration}. \ed

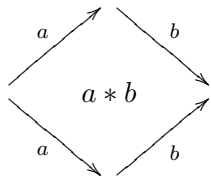
\begin{figure}
\[
\xymatrix{
&  \ar@{->}[rd]^{b}&\\
\ar@{->}[ru]^{a}\ar@{->}[rd]_{a} & a*b& \\
&\ar@{->}[ru]_{b}&}
\]
\caption{Sequential execution of $a$ and $b$}
\label{seqab}
\end{figure}

\begin{figure}
\[
\xymatrix{
&  \ar@{->}[rd]^{b}&\\
\ar@{->}[ru]^{a}\ar@{->}[rd]_{b} & a*b=b*a & \\
&\ar@{->}[ru]_{a}&}
\]
\caption{Concurrent execution of $a$ and $b$}
\label{concab2}
\end{figure}
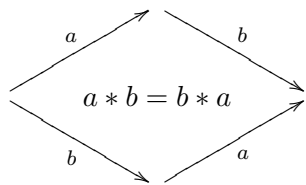

Figure~\ref{seqab} represents a labelled flow such that each execution
path from the initial state to the final state corresponds to the
\emph{sequential} execution of $a$ and $b$. The continuous deformation
between the top execution path and the bottom execution path means
that it is impossible to observe what execution path is really chosen.

Figure~\ref{concab2} represents a labelled flow corresponding to the
\emph{concurrent} execution of $a$ and $b$. Note that we need the
algebraic relation $a*b=b*a$. Thus, one must have $\sigma(a,0)=a$ and
$\sigma(b,0)=b$ by definition of $?^\sigma\Sigma$.

\section{Geometric realization of a precubical set}
\label{rea}

\begin{figure}
\[
\xymatrix{
(\widehat{0},\widehat{0}) \fr{(\widehat{0},*)}\fd{(*,\widehat{0})}\ar@{->}[ddrr]^-{(*,*)} && (\widehat{0},\widehat{1})\fd{(*,\widehat{1})}\\
&& \\
(\widehat{1},\widehat{0}) \fr{(\widehat{1},*)}&& (\widehat{1},\widehat{1})}
\]
\caption{The flow $|\square[2]|_{bad} = \{\widehat{0}<\widehat{1}\}^2$}
\label{cube2}
\end{figure}
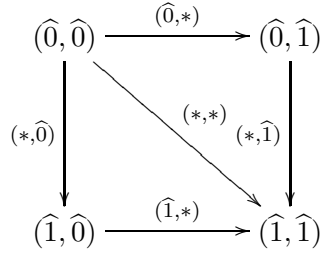

A state of the flow associated with the poset
$\{\widehat{0}<\widehat{1}\}^n$ (i.e. the product of $n$ copies of
$\{\widehat{0}<\widehat{1}\}$) is denoted by an $n$-uple of elements
of $\{\widehat{0},\widehat{1}\}$. By convention,
$\{\widehat{0}<\widehat{1}\}^0=\{()\}$. The unique morphism/execution
path from $(x_1,\dots,x_n)$ to $(y_1,\dots,y_n)$ is denoted by an
$n$-uple $(z_1,\dots,z_n)$ of $\{\widehat{0},\widehat{1},*\}$ with
$z_i=x_i$ if $x_i=y_i$ and $z_i=*$ if $x_i<y_i$. For example in the
flow $\{\widehat{0}<\widehat{1}\}^2$ (cf. Figure~\ref{cube2}), one has
the algebraic relation $(*,*) = (\widehat{0},*)*(*,\widehat{1}) =
(*,\widehat{0}) * (\widehat{1},*)$.

Let $\square \rightarrow \poset \subset \dtop$ be the functor defined
on objects by the mapping $[n]\mapsto \{\widehat{0}<\widehat{1}\}^n$
and on morphisms by the mapping \[\delta_i^\alpha \mapsto \lp
(\epsilon_1, \dots, \epsilon_{n-1}) \mapsto (\epsilon_1, \dots,
\epsilon_{i-1}, \alpha, \epsilon_i, \dots, \epsilon_{n-1})\rp,\] where
the $\epsilon_i$'s are elements of $\{\widehat{0},\widehat{1},*\}$.

Note by Yoneda's lemma, that one has the natural isomorphism of
precubical sets \[K\iso \liminj_{\square[n]\rightarrow K}
\square[n].\] The functor $[n] \mapsto \{\widehat{0}<\widehat{1}\}^n$
from $\square$ to $\dtop$ induces a realization functor from
$\square^{op}\set$ to $\dtop$ defined
by \[\boxed{|K|_{bad}:=\liminj_{\square[n]\rightarrow K}
  \{\widehat{0}<\widehat{1}\}^n}\] However, this functor is \emph{not}
the good realization functor from precubical sets to flows.  Indeed,
the flow $|\square[2]|_{bad} \iso \{\widehat{0} < \widehat{1}\} \p
\{\widehat{0} < \widehat{1}\}$ is isomorphic to the flow of
Figure~\ref{cube2}. It is then not difficult to see that there is
exactly one non-constant execution path from the initial state to the
final state of $|\de\square[3]|_{bad}$.  More generally, one has:

\bth \label{pbpb} For $n\geq 3$, the canonical inclusion
$\de\square[n] \rightarrow \square[n]$ induces the isomorphism of
flows $|\de\square[n]|_{bad} \iso |\square[n]|_{bad}$.  \eth

The flow $|\de\square[2]|_{bad}$ is obtained from the flow
$|\square[2]|_{bad}$ by removing the algebraic relation
$(\widehat{0},*)*(*,\widehat{1}) = (*,\widehat{0}) * (\widehat{1},*)$
(cf. Figure~\ref{sqsq}). Therefore the flows $|\de\square[2]|_{bad}$
and $|\square[2]|_{bad}$ are not isomorphic.

\begin{figure}
\[
\xymatrix{
(\widehat{0},\widehat{0}) \fr{(\widehat{0},*)}\fd{(*,\widehat{0})}  && (\widehat{0},\widehat{1})\fd{(*,\widehat{1})}\\
&(\widehat{0},*)*(*,\widehat{1}) \neq (*,\widehat{0}) *
  (\widehat{1},*)& \\
(\widehat{1},\widehat{0}) \fr{(\widehat{1},*)}&& (\widehat{1},\widehat{1})}
\]
\caption{The flow $|\de\square[2]|_{bad}$}
\label{sqsq}
\end{figure}

\bpf Let $n \geq 3$. One has $|\square[n]|_{bad} \iso \{\widehat{0} <
\widehat{1}\}^n$ since the comma category $\square\downarrow
\square[n]$ has a terminal object $\square[n] \rightarrow \square[n]$.
The map
\[
|\de\square[n]|_{bad} =
\int^{[r]}(\de\square[n])_r.\{\widehat{0} < \widehat{1}\}^r
\rightarrow |\square[n]|_{bad}
\]
then induces a commutative diagram of flows
\[
\bigsqcup_{[q]\rightarrow [p]}
(\de\square[n])_p.\{\widehat{0}<\widehat{1}\}^q \rightrightarrows
\bigsqcup _{[r]} (\de\square[n])_r.\{\widehat{0}<\widehat{1}\}^r
\rightarrow \{\widehat{0} < \widehat{1}\}^n.\] Consider a commutative
diagram of flows
\[
\bigsqcup_{[q]\rightarrow [p]}
(\de\square[n])_p.\{\widehat{0}<\widehat{1}\}^q \rightrightarrows
\bigsqcup _{[r]} (\de\square[n])_r.\{\widehat{0}<\widehat{1}\}^r
\rightarrow Z.\] It suffices to prove that one then has a unique
factorization
\[
\bigsqcup_{[q]\rightarrow [p]}
(\de\square[n])_p.\{\widehat{0}<\widehat{1}\}^q \rightrightarrows
\bigsqcup _{[r]} (\de\square[n])_r.\{\widehat{0}<\widehat{1}\}^r
\rightarrow \{\widehat{0}<\widehat{1}\}^n \rightarrow Z\] to complete
the proof of the proposition. By definition of the precubical set $\de
\square[n]$, one has the commutative diagram
\[
\xymatrix{ \bigsqcup\limits_{[q]\rightarrow [p],p<n}
  \square([p],[n]).\{\widehat{0}<\widehat{1}\}^q \rightrightarrows
  \bigsqcup\limits_{[r],r<n}
  \square([r],[n]).\{\widehat{0}<\widehat{1}\}^r
  \ar@{->}[rrr]^-{\bigsqcup\limits_{f\in \square([r],[n]),0\leq
      r<n}\phi_f}&&& Z,}\] where $\phi_f\colon
\{\widehat{0}<\widehat{1}\}^r \rightarrow Z$ is a morphism of flows
from the copy of $\{\widehat{0}<\widehat{1}\}^r$ indexed by $f\in
\square([r],[n])$ to $Z$. The commutativity of the above diagram means
that for any $(\epsilon_1,\dots,\epsilon_q) \in
\{\widehat{0}<\widehat{1}\}^q$, one has
\begin{equation} \label{commutatif} \boxed{\phi_{f\circ
      g}(\epsilon_1,\dots,\epsilon_q) =
\phi_f(|g|_{bad}(\epsilon_1,\dots,\epsilon_q))}
\end{equation}
where $g\colon [q] \rightarrow [p]$ and $f\colon [p] \rightarrow [n]$
are two morphisms of the category $\square$ with $p<n$. Suppose that
this factorization $h\colon \{\widehat{0}<\widehat{1}\}^n \rightarrow
Z$ exists.  Then necessarily, \beas h(\epsilon_1,\dots,\epsilon_n) &=&
h(s(\epsilon_1),\epsilon_2,\dots,\epsilon_n)*h(\epsilon_1,t(\epsilon_2),\dots,t(\epsilon_n))\\
& =&
\phi_{\delta^{s(\epsilon_1)}_1}(\epsilon_2,\dots,\epsilon_n)*\phi_{\delta_n^{t(\epsilon_n)}\dots\delta_2^{t(\epsilon_2)}}(\epsilon_1).
\eeas So there is at most one such map $h$.  Let
\[h(\epsilon_1,\dots,\epsilon_n)=\phi_{\delta^{s(\epsilon_1)}_1}(\epsilon_2,\dots,\epsilon_n)*\phi_{\delta_n^{t(\epsilon_n)}\dots\delta_2^{t(\epsilon_2)}}(\epsilon_1).\]
One has
$s(h(\epsilon_1,\dots,\epsilon_n))=s(\phi_{\delta^{s(\epsilon_1)}_1}(\epsilon_2,\dots,\epsilon_n))=\phi_{\delta^{s(\epsilon_1)}_1}(s(\epsilon_2),\dots,s(\epsilon_n))$
since $\phi_{\delta^{s(\epsilon_1)}_1}$ is a morphism of flows and
\beas
h(s(\epsilon_1),\dots,s(\epsilon_n)) &=& \phi_{\delta^{s(\epsilon_1)}_1}(s(\epsilon_2),\dots,s(\epsilon_n))*\phi_{\delta_n^{s(\epsilon_n)}\dots\delta_2^{s(\epsilon_2)}}(s(\epsilon_1))\\
&=&  \phi_{\delta^{s(\epsilon_1)}_1}(s(\epsilon_2),\dots,s(\epsilon_n)) * \phi_{\delta_n^{s(\epsilon_n)}\dots\delta_3^{s(\epsilon_3)}}(s(\epsilon_1),s(\epsilon_2)) \\
&=& \phi_{\delta^{s(\epsilon_1)}_1}(s(\epsilon_2),\dots,s(\epsilon_n)) * \phi_{\delta_n^{s(\epsilon_n)}\dots\delta_3^{s(\epsilon_3)}\delta_1^{s(\epsilon_1)}}(s(\epsilon_2))\\
&=& \phi_{\delta^{s(\epsilon_1)}_1}(s(\epsilon_2),\dots,s(\epsilon_n)) * \phi_{\delta^{s(\epsilon_1)}_1}(s(\epsilon_2),\dots,s(\epsilon_n)) \\
&=&
\phi_{\delta^{s(\epsilon_1)}_1}(s(\epsilon_2),\dots,s(\epsilon_n)).
\eeas Thus,
$s(h(\epsilon_1,\dots,\epsilon_n))=h(s(\epsilon_1),\dots,s(\epsilon_n))$
and, in the same way, one can show the equality
$t(h(\epsilon_1,\dots,\epsilon_n))=h(t(\epsilon_1),\dots,t(\epsilon_n))$. It
remains to prove
that \[h(\epsilon_1*\epsilon'_1,\dots,\epsilon_n*\epsilon'_n)=h(\epsilon_1,\dots,\epsilon_n)
* h(\epsilon'_1,\dots,\epsilon'_n).\] It is equivalent to proving the
equality
\begin{multline*}
\phi_{\delta^{s(\epsilon_1)}_1}(\epsilon_2,\dots,\epsilon_n) *
\phi_{\delta^{s(\epsilon_1)}_1}(\epsilon'_2,\dots,\epsilon'_n) *
\phi_{\delta_n^{t(\epsilon'_n)}\dots\delta_2^{t(\epsilon'_2)}}(\epsilon_1) *
\phi_{\delta_n^{t(\epsilon'_n)}\dots\delta_2^{t(\epsilon'_2)}}(\epsilon'_1) \\ 
=
\phi_{\delta^{s(\epsilon_1)}_1}(\epsilon_2,\dots,\epsilon_n)*\phi_{\delta_n^{t(\epsilon_n)}\dots\delta_2^{t(\epsilon_2)}}(\epsilon_1)
*
\phi_{\delta^{s(\epsilon'_1)}_1}(\epsilon'_2,\dots,\epsilon'_n)*\phi_{\delta_n^{t(\epsilon'_n)}\dots\delta_2^{t(\epsilon'_2)}}(\epsilon'_1).
\end{multline*}
It therefore remains to prove the equality
\[\phi_{\delta^{s(\epsilon_1)}_1}(\epsilon'_2,\dots,\epsilon'_n) * 
\phi_{\delta_n^{t(\epsilon'_n)}\dots\delta_2^{t(\epsilon'_2)}}(\epsilon_1) =
\phi_{\delta_n^{t(\epsilon_n)}\dots\delta_2^{t(\epsilon_2)}}(\epsilon_1) *
\phi_{\delta^{s(\epsilon'_1)}_1}(\epsilon'_2,\dots,\epsilon'_n)\]
or the equality 
\[\phi_{\delta^{s(\epsilon_1)}_1}(\epsilon'_2,\dots,\epsilon'_n) * 
\phi_{\delta_n^{t(\epsilon'_n)}\dots\delta_2^{t(\epsilon'_2)}}(\epsilon_1) =
\phi_{\delta_n^{s(\epsilon'_n)}\dots\delta_2^{s(\epsilon'_2)}}(\epsilon_1) *
\phi_{\delta^{t(\epsilon_1)}_1}(\epsilon'_2,\dots,\epsilon'_n)\] since $s(\epsilon'_i)=t(\epsilon_i)$ for
all $1 \leq i \leq n$.  One then has{\footnotesize
\begin{align*}
& \phi_{\delta^{s(\epsilon_1)}_1}(\epsilon'_2,\dots,\epsilon'_n) * 
\phi_{\delta_n^{t(\epsilon'_n)}\dots\delta_2^{t(\epsilon'_2)}}(\epsilon_1) & \\
& = \phi_{\delta^{s(\epsilon_1)}_1}(s(\epsilon'_2),\epsilon'_3,\dots,\epsilon'_n) * \phi_{\delta^{s(\epsilon_1)}_1}(\epsilon'_2,t(\epsilon'_3),\dots,t(\epsilon'_n)) * 
\phi_{\delta_n^{t(\epsilon'_n)}\dots\delta_2^{t(\epsilon'_2)}}(\epsilon_1) & 
\\
& = \phi_{\delta^{s(\epsilon'_2)}_2}(s(\epsilon_1),\epsilon'_3,\dots,\epsilon'_n) *  \phi_{\delta_n^{t(\epsilon'_n)}\dots\delta_3^{t(\epsilon'_3)}}(s(\epsilon_1),\epsilon'_2)* \phi_{\delta_n^{t(\epsilon'_n)}\dots\delta_3^{t(\epsilon'_3)}}(\epsilon_1,t(\epsilon'_2)) & 
\hbox{ by (\ref{commutatif})}
\\
& = \phi_{\delta^{s(\epsilon'_2)}_2}(s(\epsilon_1),\epsilon'_3,\dots,\epsilon'_n) *  \phi_{\delta_n^{t(\epsilon'_n)}\dots\delta_3^{t(\epsilon'_3)}}(\epsilon_1,s(\epsilon'_2))* \phi_{\delta_n^{t(\epsilon'_n)}\dots\delta_3^{t(\epsilon'_3)}}(t(\epsilon_1),\epsilon'_2) & 
\\
& = \phi_{\delta^{s(\epsilon'_2)}_2}(s(\epsilon_1),\epsilon'_3,\dots,\epsilon'_n) * \phi_{\delta^{s(\epsilon'_2)}_2}(\epsilon_1,t(\epsilon'_3),\dots,t(\epsilon'_n)) *  \phi_{\delta_n^{t(\epsilon'_n)}\dots\delta_3^{t(\epsilon'_3)}}(t(\epsilon_1),\epsilon'_2)& 
\hbox{ by (\ref{commutatif})}
\\
& = \phi_{\delta^{s(\epsilon'_2)}_2}(\epsilon_1,s(\epsilon'_3),\dots,s(\epsilon'_n)) * \phi_{\delta^{s(\epsilon'_2)}_2}(t(\epsilon_1),\epsilon'_3,\dots,\epsilon'_n) *  \phi_{\delta_1^{t(\epsilon_1)}}(\epsilon'_2,t(\epsilon'_3),\dots,t(\epsilon'_n)) & 
\\ 
& = \phi_{\delta^{s(\epsilon'_n)}_n\dots\delta^{s(\epsilon'_2)}_2}(\epsilon_1) * \phi_{\delta^{t(\epsilon_1)}_1}(s(\epsilon'_2),\epsilon'_3,\dots,\epsilon'_n) *  \phi_{\delta_1^{t(\epsilon_1)}}(\epsilon'_2,t(\epsilon'_3),\dots,t(\epsilon'_n)) & 
\hbox{ by (\ref{commutatif})}
\\
& =
\phi_{\delta^{s(\epsilon'_n)}_n\dots\delta^{s(\epsilon'_2)}_2}(\epsilon_1)
*
\phi_{\delta^{t(\epsilon_1)}_1}(\epsilon'_2,\epsilon'_3,\dots,\epsilon'_n).  & 
\end{align*}}\epf

The origin of the problem is the algebraic relation
$(\widehat{0},*)*(*,\widehat{1}) = (*,\widehat{0}) * (\widehat{1},*)$
in $\{\widehat{0}<\widehat{1}\}^2$. To get a correct behaviour, it is
necessary to replace algebraic relations like
$(\widehat{0},*)*(*,\widehat{1}) = (*,\widehat{0}) * (\widehat{1},*)$
by homotopies between $(\widehat{0},*)*(*,\widehat{1})$ and
$(*,\widehat{0}) * (\widehat{1},*)$. The good definition of the
realization of a precubical set as flow therefore requires the weak
S-homotopy model structure constructed in~\cite{model3}.

\bd Let $K$ be a precubical set. By definition, the \emph{geometric
  realization} of $K$ is the flow
\[\boxed{|K| := \liminj_{\square[n]\rightarrow
    K} (\{\widehat{0}<\widehat{1}\}^n)^{\textit{cof}}} \]
\ed

\bp \label{samestates} Let $K$ be a precubical set. Then one has the
natural bijections of sets $K_0 \iso |K|_{bad}^0 \iso |K|^0$. \ep

\bpf The functor $X\mapsto X^0$ from $\dtop$ to $\set$ is
colimit-preserving since it is a left adjoint: one has the natural
bijection $\set(X^0,S)\iso \dtop(X,\widehat{S})$, where $\widehat{S}$
is the flow defined by $\widehat{S}^0 = S$ and
$\P_{\alpha,\beta}\widehat{S} = \{(\alpha,\beta)\}$ for every
$(\alpha,\beta)\in S\p S$. So \[|K|_{bad}^0 \iso |K|^0 \iso
\liminj_{\square[n]\rightarrow K} (\{\widehat{0}<\widehat{1}\}^n)^0
\iso \liminj_{\square[n]\rightarrow K} \square[n]_0 =
K_0.\] \epf

The following propositions help to understand the difference between
the bad geometric realization and the good geometric realization
functors:

\bp \label{decompcof} For every $n\geq 0$, the map of flows $|\de
\square[n]| \rightarrow |\square[n]|$ is a cofibration.  \ep

Note that it is also true that the map of flows $|\de \square[n]|_{bad}
\rightarrow |\square[n]|_{bad}$ is a co\-fi\-bration for all $n\geq 3$
since an isomorphism is a cofibration. However, the map $|\de
\square[2]|_{bad} \rightarrow |\square[2]|_{bad}$ is not a
cofibration. 

\bpf One has $|\square[n]|
=(\{\widehat{0}<\widehat{1}\}^n)^{\textit{cof}}$ by definition of the
realization functor for every $n\geq 1$.  Therefore the flow
$|\square[n]|$ is a $\{\glob(\mathbf{S}^{p-1}) \subset
\glob(\mathbf{D}^{p}), p\geq 0\}$-cell complex. The flow
$|\de\square[n]|$ is a $\{\glob(\mathbf{S}^{p-1}) \subset
\glob(\mathbf{D}^{p}), p\geq 0\}$-cell subcomplex, the one obtained by
removing the cells $\glob(\mathbf{S}^{p-1}) \subset
\glob(\mathbf{D}^{p})$ of $|\square[n]|$ whose attaching map sends the
initial and final states of $\glob(\mathbf{S}^{p-1})$ to respectively
the initial and final states of $|\square[n]|$.  Thus, the map $|\de
\square[n]| \rightarrow |\square[n]|$ is a cofibration.  \epf

\bp\label{reacof} For every precubical set $K$, the geometric
realization $|K|$ is a cofibrant flow. \ep

\bpf Let $K$ be a precubical set. The map $\varnothing \rightarrow
K_0$ is a transfinite composition of pushouts of the map
$\de\square[0]=\varnothing \rightarrow \square[0]=\{0\}$. The passage
from $K_{\leq n-1}$ to $K_{\leq n}$ for $n\geq 1$ is done by the
following pushout diagram:
\[
\xymatrix{
\bigsqcup_{x\in K_n} \de\square[n] \fr{} \fd{} && K_{\leq n-1} \fd{} \\
&& \\
\bigsqcup_{x\in K_n} \square[n] \fr{} && K_{\leq n},\cocartesien}
\] 
where the map $\de\square[n] \rightarrow K_{\leq n-1}$ indexed by
$x\in K_n$ is induced by the $(n-1)$-shell $\de\square[n] \subset
\square[n] \stackrel{x}\rightarrow K$. Therefore $K$ is a
$\{\de\square[n] \subset \square[n],n\geq 0\}$-cell complex.  The
proposition is then a consequence of Proposition~\ref{decompcof}. \epf

Possible references for the classifying space of a small category
are~\cite{MR0338129, MR0232393}. For homotopy theory of posets,
see~\cite{poset_tool}. If $P$ is a poset, then $\Delta(P)$ denotes the
\emph{order complex} associated with $P$. Recall that the order
complex is a simplicial complex having $P$ as underlying set and
having the subsets $\{x_0,x_1,\dots,x_n\}$ with $x_0<x_1<\dots<x_n$ as
$n$-simplices~\cite{MR493916}.  Such a simplex will be denoted by
$(x_0,x_1,\dots,x_n)$. The order complex $\Delta(P)$ can be viewed as
a poset ordered by the inclusion, and therefore as a small category.
The corresponding category will be denoted in the same way. The
opposite category $\Delta(P)^{op}$ is freely generated by the
morphisms $\de_i\colon (x_0,\dots,x_n) \longrightarrow
(x_0,\dots,\widehat{x_i},\dots,x_n)$ for $0\leq i\leq n$ and by the
simplicial relations $\de_i\de_j=\de_{j-1}\de_i$ for any $i<j$, where
the notation $\widehat{x_i}$ means that $x_i$ is removed.  The order
complex $\Delta(P)$ corresponds to the barycentric subdivision of $P$.
Hence $P$ and $\Delta(P)$ have the same homotopy type.

\bth \label{nopb} Let $n\geq 1$. The topological space
$\P_{\widehat{0}\dots\widehat{0},\widehat{1}\dots\widehat{1}}
|\de\square[n]|$ of non-constant execution paths from the initial
state of $|\de\square[n]|$ to its final state is homotopy equivalent
to $\mathbf{S}^{n-2}$. \eth

\bpf The statement is obvious for $n=1,2$. Let $n\geq 3$. Following
the notation of~\cite[\S7]{4eme}, let us denote by
$\Delta^{ext}(\{\widehat{0}<\widehat{1}\}^n)$ the full subcategory of
$\Delta(\{\widehat{0}<\widehat{1}\}^n)$ consisting of the simplices
$(\alpha_0,\dots,\alpha_p)$ such that $\widehat{0}\dots\widehat{0} =
\alpha_0$ and $\widehat{1}\dots\widehat{1} = \alpha_p$.  The simplex
$(\widehat{0}\dots\widehat{0},\widehat{1}\dots\widehat{1})$ is a
terminal object of $\Delta^{ext}(\{\widehat{0}<\widehat{1}\}^n)^{op}$.
This category can be equipped with a structure of direct Reedy
category by~\cite[Proposition 7.3]{4eme}. Consider the diagram of
topological spaces \[\mathcal{F}|\de\square[n]|:
\Delta^{ext}(\{\widehat{0}<\widehat{1}\}^n)^{op}\backslash
\{(\widehat{0}\dots\widehat{0},\widehat{1}\dots\widehat{1})\}
\longrightarrow \top\] defined by
\[\mathcal{F}|\de\square[n]|(\alpha_0,\dots,\alpha_p) = \P_{\alpha_0,\alpha_1}
|\de\square[n]| \p \dots \p \P_{\alpha_{p-1},\alpha_p}
|\de\square[n]|\] on objects and by the composition law of the flow
$|\de\square[n]|$ on arrows (cf.~\cite[Proposition 7.5]{4eme}. Since
the flow $|\de\square[n]|$ is cofibrant by Proposition~\ref{reacof},
the diagram of topological spaces $\mathcal{F}|\de\square[n]|$ is then
Reedy cofibrant by~\cite[Theorem 7.14]{4eme} and~\cite[Theorem
8.1]{4eme}.

The point is that one has the isomorphism of topological spaces
\[\liminj
\mathcal{F}|\de\square[n]| \iso
\P_{\widehat{0}\dots\widehat{0},\widehat{1}\dots\widehat{1}}
|\de\square[n]|\] since the latter topological space is freely
generated by the non-constant execution paths of each face of the
boundary.  This colimit is actually a homotopy colimit
by~\cite[Theorem 15.10.8]{ref_model2} since the category
$\Delta^{ext}(\{\widehat{0}<\widehat{1}\}^n)^{op}\backslash
\{(\widehat{0}\dots\widehat{0},\widehat{1}\dots\widehat{1})\}$ has
fibrant constants by~\cite[Definition 15.10.1]{ref_model2}
and~\cite[Proposition 15.10.2]{ref_model2}. So one has a weak homotopy
equivalence \[\holiminj \mathcal{F}|\de\square[n]|
\stackrel{\simeq}\rightarrow
\P_{\widehat{0}\dots\widehat{0},\widehat{1}\dots\widehat{1}}
|\de\square[n]|.\]

Since all topological spaces
$\mathcal{F}|\de\square[n]|(\alpha_0,\dots,\alpha_p)$ are
contractible, there is an ob\-ject\-wise weak homotopy equivalence
\[\mathcal{F}|\de\square[n]| \stackrel{\simeq}\longrightarrow
\mathbf{1}_{\Delta^{ext}(\{\widehat{0}<\widehat{1}\}^n)^{op}\backslash
  \{(\widehat{0}\dots\widehat{0},\widehat{1}\dots\widehat{1})\}},\]
where the diagram
$\mathbf{1}_{\Delta^{ext}(\{\widehat{0}<\widehat{1}\}^n)^{op}\backslash
  \{(\widehat{0}\dots\widehat{0},\widehat{1}\dots\widehat{1})\}}$ is
the terminal diagram over the small category
$\Delta^{ext}(\{\widehat{0}<\widehat{1}\}^n)^{op}\backslash
\{(\widehat{0}\dots\widehat{0},\widehat{1}\dots\widehat{1})\}$.
By~\cite[Proposition 18.1.6]{ref_model2}, the topological space
$\liminj \mathcal{F}|\de\square[n]|$ is therefore weakly homotopy
equivalent to the classifying space of the index category
$\Delta^{ext}(\{\widehat{0}<\widehat{1}\}^n)^{op}\backslash
\{(\widehat{0}\dots\widehat{0},\widehat{1}\dots\widehat{1})\}$. The
latter category is the category of simplices of the poset
$\{\widehat{0}<\widehat{1}\}^n \backslash
\{\widehat{0}\dots\widehat{0},\widehat{1}\dots\widehat{1}\}$, that is
its order complex, or topologically its barycentric subdivision.

The poset associated with the $n$-cube is isomorphic to the poset of
lattices of the faces of the $(n-1)$-simplex. More precisely, an
element of the $n$-cube is represented by a word of length $n$ in the
alphabet $\{\widehat{0},\widehat{1}\}$.  Let $\epsilon_1\dots
\epsilon_n$. Then let $\phi(\epsilon_1\dots \epsilon_n)$ be the subset
of $\{1,\dots,n-1\}$ defined by $i\in \phi(\epsilon_1\dots
\epsilon_n)$ if and only if $\epsilon_i=\widehat{1}$. So $\phi$
induces an isomorphism between the poset associated with the $n$-cube
and the order complex of $\{1<2<\dots<n-1\}$. Thus the poset
$\{\widehat{0}<\widehat{1}\}^n \backslash
\{\widehat{0}\dots\widehat{0},\widehat{1}\dots\widehat{1}\}$ has the
homotopy type of the boundary of the $(n-1)$-simplex (\cite[Example
1.1.1]{poset_tool}). So the topological spaces
$\P_{\widehat{0}\dots\widehat{0},\widehat{1}\dots\widehat{1}}
|\de\square[n]|$ and $\mathbf{S}^{n-2}$ are weakly homotopy
equivalent.  Hence we have the result since the topological space
$\P_{\widehat{0}\dots\widehat{0},\widehat{1}\dots\widehat{1}}
|\de\square[n]|$ is cofibrant. \epf

\begin{corollary} Let $n\geq 0$. The cofibration $|\de\square[n]|
  \rightarrow |\square[n]|$ is non-trivial. \end{corollary}

\bpf Corollary of Proposition~\ref{decompcof} and Theorem~\ref{nopb}
for $n\geq 2$. The proof is clear for $n=0,1$.  \epf

The underlying homotopy type of a flow is, morally speaking, the
underlying state space of a time flow, unique only up to homotopy: it
is defined in~\cite{model2} and studied in~\cite{4eme}.

\begin{corollary} \label{correct} Let $n\geq 1$. The underlying
  homotopy type of the flow $|\de\square[n]|$ is
  $\mathbf{S}^{n-1}$. \end{corollary}

Notice that the underlying homotopy type of $|\de\square[n]|_{bad}$
for $n\geq 3$ is that of $|\square[n]|_{bad}$ by Theorem~\ref{pbpb},
that is a point!

\bpf By Theorem~\ref{nopb} and~\cite[Corollary 8.7]{4eme}, the
underlying homotopy type of the flow $|\de\square[n]|$ is the
underlying homotopy type of the flow $\glob(\mathbf{S}^{n-2})$, and
the non-pointed suspension of $\mathbf{S}^{n-2}$ is precisely
$\mathbf{S}^{n-1}$.  \epf

The natural trivial fibrations $(\{\widehat{0} <
\widehat{1}\}^n)^{\textit{cof}} \longrightarrow \{\widehat{0} <
\widehat{1}\}^n$ for $n\geq 0$ induce a natural map $|K|
\longrightarrow |K|_{bad}$ for any precubical set $K$. In general, the
latter map is not a weak S-homotopy equivalence. Indeed, the
topological space
$\P_{(\widehat{0},\widehat{0},\widehat{0}),(\widehat{1},\widehat{1},\widehat{1})}
|\de\square[3]|$ is homotopy equivalent to the $1$-dimensional sphere
$\mathbf{S}^1$ by Theorem~\ref{nopb}, and it is already known by
Theorem~\ref{pbpb} that the topological space
$\P_{(\widehat{0},\widehat{0},\widehat{0}),(\widehat{1},\widehat{1},\widehat{1})}
|\de\square[3]|_{bad}$ is a singleton. However, one has:

\bp \label{flowdis} For every precubical set $K$, the flow $|K|_{bad}$
has a discrete space of non-constant execution paths; i.e.,
$\P(|K|_{bad})$ is discrete. So the natural map $|K| \rightarrow
|K|_{bad}$ is a fibration of flows.  \ep

\bpf By definition, one has $|K|_{bad}=\liminj_{\square[n]\rightarrow
  K} \{\widehat{0} < \widehat{1}\}^n$. The path space
\[
\P(\{\widehat{0} < \widehat{1}\}^n)
\]
is discrete by construction. The path space $\P(|K|_{bad})$ is
generated by the free compositions of elements of
$\liminj_{\square[n]\rightarrow K} \P(\{\widehat{0} <
\widehat{1}\}^n)$. Hence the result since a colimit of discrete spaces
is discrete.  Since every topological space is fibrant, one deduces
that the natural map $|K| \rightarrow |K|_{bad}$ is a fibration of
flows.  \epf

\section{Associating a $\sigma$-labelled precubical set with a $\sigma$-labelled flow}
\label{reaet}

\bp \label{badsigma} One has the isomorphism of flows
$|!^\sigma\Sigma|_{bad} \iso\ ?^\sigma\Sigma$. \ep

\bpf The functor $K\mapsto |K|_{bad}$ is a left adjoint since the
right adjoint is the functor $X\mapsto ([n] \mapsto
\dtop(\{\widehat{0},\widehat{1}\}^n,X))$; therefore it preserves all
colimits. So by Theorem~\ref{pbpb}, the canonical inclusion
$(!^\sigma\Sigma)_{\leq 2} \subset !^\sigma\Sigma$, which is a relative
$\{\de\square[n]\subset \square[n],n\geq 3\}$-cell complex, induces an
isomorphism of flows $|(!^\sigma\Sigma)_{\leq 2}|_{bad} \iso
|!^\sigma\Sigma|_{bad}$. By Proposition~\ref{samestates}, the set
$|!^\sigma\Sigma|_{bad}^0$ is a singleton.  Each element of $\Sigma$
generates a non-constant execution path of $|!^\sigma\Sigma|_{bad}$
and each algebraic relation $a*b=b*a$ corresponds to a $2$-cube
$(a,b)$ of $!^\sigma\Sigma$ (with $a\leq b$).
So one has the
isomorphism $|(!^\sigma\Sigma)_{\leq 2}|_{bad} \iso ?^\sigma\Sigma$.
\epf

Let $K\rightarrow !^\sigma\Sigma$ be a (resp.\ decorated)
$\sigma$-labelled precubical set. Then the composition $|K|
\rightarrow |!^\sigma\Sigma| \rightarrow |!^\sigma\Sigma|_{bad} \iso
?^\sigma\Sigma$ gives rise to a (resp.\ decorated) $\sigma$-labelled
flow. Together with Section~\ref{semccsprecube}, one obtains an
interpretation $|\square|[P|]|$ for every process name $P$ in terms of
$\sigma$-labelled flows. The flow $|\square|[P|]|$ is cofibrant for
every process name $P$ of $\proc_\Sigma$ by Proposition~\ref{reacof}.

\appendix

\section{Associativity of the synchronized tensor product}
\label{associativity}

For the sequel, the category of small categories is denoted by
$\cat$. Let $H\colon I \longrightarrow \cat$ be a functor from a small
category $I$ to $\cat$. ``The'' \emph{Grothendieck construction} $I
\intop H$ is the category defined as follows~\cite{MR510404}: the
objects are the pairs $(i,a)$ where $i$ is an object of $I$ and $a$ is
an object of $H(i)$; a morphism $(i,a) \rightarrow (j,b)$ consists in
a map $\phi\colon i \rightarrow j$ and in a map $h \colon H(\phi)(a)
\rightarrow b$.

\begin{lemma} \label{colim_id} Let $I$ be a small category, and 
  $i\mapsto K^i$ be a functor from $I$ to the category of
  $\sigma$-labelled precubical sets.  Let $K = \liminj_i K^i$. Let
  $H\colon I \rightarrow \cat$ be the func\-tor defined by $H(i) = \square
  \downarrow K^i$.  Then the functor $\iota \colon  I\intop H \rightarrow
  \square \downarrow K$ defined by $\iota(i,\square[m] \rightarrow
  K^i) = (\square[m] \rightarrow K)$ is final in the sense 
  of~\cite{MR1712872}; that is to say the comma category $k \downarrow
  \iota$ is nonempty and connected for all objects $k$ of $\square
  \downarrow K$.  \end{lemma}

\bpf Let $k\colon \square[m] \rightarrow K$ be an object of $\square
\downarrow K$. Then an object of the comma category $k \downarrow
\iota$ consists of a pair $((j,\square[n] \rightarrow K^j), \phi)$,
where $\phi\colon \square[m] \rightarrow \square[n]$ is a map of
precubical sets such that the following diagram is commutative:
\[
\xymatrix{
\square[m]\fr{k} \fd{\phi}&& K \\
&& \\
\square[n] \fr{} && K^j. \fu{}}
\] 
By Yoneda, $k \in K_m = \liminj_i (K^i)_m$. So there exists $i$ such
that $k$ factors as a composite $\square[m] \rightarrow K^i
\rightarrow K$, and the pair $((i,\square[m] \rightarrow K^i),
\id_{\square[m]})$ is an object of $k \downarrow \iota$. Thus the
latter category is not empty. One has to prove that the latter
category is connected. Let $((j,\square[n] \rightarrow K^j), \phi)$ be
an object of $k \downarrow \iota$. The map $\phi$ induces a map from
$((j,\square[m] \rightarrow K^j), \id_{\square[m]})$ to
$((j,\square[n] \rightarrow K^j), \phi)$. Thus one can suppose $m=n$
and $\phi = \id_{\square[m]}$ without loss of generality. Let
$((j_1,\square[m] \rightarrow K^{j_1}), \id_{\square[m]})$ and
$((j_2,\square[m] \rightarrow K^{j_2}), \id_{\square[m]})$ be two
objects of $k \downarrow \iota$. By definition, the two composites
$\square[m] \stackrel{k_1}\rightarrow K^{j_1} \rightarrow K$ and
$\square[m] \stackrel{k_2}\rightarrow K^{j_2} \rightarrow K$ are equal
to $k$. Since $k_1 \in (K^{j_1})_m$ and $k_2 \in (K^{j_2})_m$ have the
same images $k\in K_m$, there exists $j_3$ and two morphisms $j_1
\rightarrow j_3$ and $j_2 \rightarrow j_3$ such that the images of
$k_1$ and $k_2$ in $(K^{j_3})_m$ are equal. Therefore one has the
zig-zag of maps of $k \downarrow \iota$ {\small\[((j_1,\square[m]
  \rightarrow K^{j_1}), \id_{\square[m]}) \rightarrow ((j_3,\square[m]
  \rightarrow K^{j_3}), \id_{\square[m]}) \leftarrow ((j_2,\square[m]
  \rightarrow K^{j_2}), \id_{\square[m]}).\]} Thus the comma category
$k \downarrow \iota$ is connected.  \epf

\bp \label{la} Let $I$ be a small category. Let $i\mapsto K^i$ be a
functor from $I$ to the category of $\sigma$-labelled precubical
sets, and let $L$ be a $\sigma$-labelled precubical set.  Then one has the
natural isomorphism $(\liminj K^i) \ot_\sigma L \iso \liminj (K^i
\ot_\sigma L)$. \ep

\bpf Let $K = \liminj K^i$. By definition, one has the isomorphism 
\[ \liminj (K^i \ot_\sigma L) \iso \liminj_i
\liminj_{\square[m]\rightarrow K^i} \liminj_{\square[n]\rightarrow L}
\square[m] \ot_\sigma \square[n].\] Consider the 
functor $H\colon I
\longrightarrow \cat$ defined by $H(i) = \square \downarrow K^i$.
Consider the functor $F_i\colon H(i) \longrightarrow
\square^{op}\set
\downarrow !^\sigma\Sigma$ defined by
\[
F_i(\square[m] \rightarrow K^i) = \liminj_{\square[n]\rightarrow L}
\square[m] \ot_\sigma \square[n].
\]
Consider the functor $F\colon I\intop H
\longrightarrow \square^{op}\set \downarrow !^\sigma\Sigma$ defined by
\[F(i,\square[m] \rightarrow K^i) = \liminj_{\square[n]\rightarrow L}
\square[m] \ot_\sigma \square[n].\] Then the composite $H(i) \subset
I\intop H \rightarrow \square^{op}\set \downarrow !^\sigma\Sigma$ is
exactly $F_i$. Therefore one has the isomorphism 
\[\liminj_i \liminj_{\square[m]\rightarrow
  K^i} \liminj_{\square[n]\rightarrow L} \square[m] \ot_\sigma
\square[n] \iso \liminj_{(i,\square[m]\rightarrow K^i)}
\liminj_{\square[n]\rightarrow L} \square[m] \ot_\sigma \square[n]\]
by~\cite[Proposition 40.2]{monographie_hocolim}. The functor
$\iota\colon 
I\intop H \rightarrow \square \downarrow K$ defined by  
$\iota(i,\square[m] \rightarrow K^i) = (\square[m] \rightarrow K)$ is
final in the sense of~\cite{MR1712872} by Lemma~\ref{colim_id}.
Therefore by~\cite[p.\ 213, Theorem 1]{MR1712872} or~\cite[Theorem
14.2.5]{ref_model2},
one has the isomorphism
\[\liminj_{(i,\square[m]\rightarrow K^i)}
\liminj_{\square[n]\rightarrow L} \square[m] \ot_\sigma \square[n] \iso
\liminj_{\square[m]\rightarrow K} \liminj_{\square[n]\rightarrow L}
\square[m] \ot_\sigma \square[n] =: K \ot_\sigma L.\]
\epf 

\bp \label{ass0} Let $p,q,r\geq 0$. Let $\square[p]$, $\square[q]$ and
$\square[r]$ be three $\sigma$-labelled full cubes. Then one has an
isomorphism of $\sigma$-labelled precubical sets \[\square[p]
\ot_\sigma(\square[q] \ot_\sigma \square[r]) \iso (\square[p]
\ot_\sigma\square[q]) \ot_\sigma \square[r].\] \ep

\bpf[Sketch of proof] By Proposition~\ref{res1}, one has \[(\square[p]
\ot_\sigma(\square[q] \ot_\sigma \square[r]))_{\leq 1} \iso
((\square[p] \ot_\sigma\square[q]) \ot_\sigma \square[r)_{\leq 1}.\]
The two $\sigma$-labelled precubical sets \[\square[p]
\ot_\sigma(\square[q] \ot_\sigma \square[r])\] and \[(\square[p]
\ot_\sigma\square[q]) \ot_\sigma \square[r]\] have the same higher
dimensional cubes parametrized by the same non-twisted maps since the
synchronization algebra $\sigma$ is associative.  The difference with
the cases of $\square[p] \ot_\sigma\square[q]$ or $\square[q]
\ot_\sigma \square[r]$ is that a coordinate may occur three times
because, for example, an action of $\square[p]$ may synchronize with
an action of $\square[q] \ot_\sigma \square[r]$ synchronizing an
action of $\square[q]$ and an action of $\square[r]$.  \epf

As a corollary, one obtains:

\bth Let $K$, $L$ and $M$ be three $\sigma$-labelled precubical
sets. Then there exists a natural isomorphism of $\sigma$-labelled
precubical sets
\[
K \ot_\sigma (L \ot_\sigma M) \iso (K \ot_\sigma L)
\ot_\sigma M.
\]
\eth

\bpf One has 
\begin{align*} 
  & K \ot_\sigma (L \ot_\sigma M) &\\
& \iso
  \lp \liminj_{\square[p]\rightarrow K} \square[p] \rp \ot_\sigma \lp
  \liminj_{\square[q]\rightarrow L} \liminj_{\square[r]\rightarrow
    M}  \square[q] \ot_\sigma \square[r]\rp & \\
  & \iso \liminj_{\square[p]\rightarrow
    K}\liminj_{\square[q]\rightarrow L} \liminj_{\square[r]\rightarrow
    M} \lp \square[p] \ot_\sigma(\square[q] \ot_\sigma \square[r])\rp & \hbox{
    by Proposition~\ref{la}}\\
  & \iso \liminj_{\square[p]\rightarrow
    K}\liminj_{\square[q]\rightarrow L} \liminj_{\square[r]\rightarrow
    M} \lp (\square[p] \ot_\sigma\square[q]) \ot_\sigma \square[r] \rp& \hbox{
    by Proposition~\ref{ass0}}\\
  &\iso (K \ot_\sigma L) \ot_\sigma M.
\end{align*} 

\epf

\end{document}